%-----------------------------------------------
%            ENCODING ALGEBRAIC POWER SERIES    VS HERWIG FEB 24  2014
%-----------------------------------------------
\input amssym.tex
\parindent=0cm 
\hoffset=1truecm
\hsize=13truecm
\vsize=23.5truecm
\baselineskip 13pt
\long\def\ignore#1\recognize{}

%-----------------------------------------------
%    FONTNAMES  FOR  TEXSHOP
%-----------------------------------------------

%\ignore

\font\Times=ptmr at 10pt
\font\Bf=ptmb at 12pt
\font\bf=ptmb at 10pt

\font\smallTimes=ptmr at 9 pt
\font\it=ptmri at 10pt

\font\eightrm=cmr8

\Times

%\recognize
%-----------------------------------------------
%                MACROS
%-----------------------------------------------

\def\litem{\par\noindent\hangindent=\parindent\ltextindent}
\def\ltextindent#1{\hbox to \hangindent{#1\hss}\ignorespaces}
\long\def\ignore#1\recognize{}

\def\big{\bigskip}
\def\med{\medskip}
\def\ds{\displaystyle}
\def\hs{\hskip}
\def\vs{\vskip}

\def\cl{\centerline}
\def\ol{\overline}
\def\sm{\setminus}

\def\cd{\cdot}
\def\to{,\ldots,}
\def\sub{\subseteq}
\def\{{\lbrace}
\def\}{\rbrace}

\def\isom{\cong}
\def\map{\rightarrow}
\def\inv{^{-1}}
\def\N{{\Bbb N}}

\def\Z{{\Bbb Z}}
\def\A{{\Bbb A}}

\def\P{{\cal P}}
\def\abs#1{\vert#1\vert}
\def\>{\rangle}
\def\<{\langle}

\def\a{\alpha}
\def\b{\beta}

\def\d{\delta}
\def\e{\varepsilon}
\def\eps{\varepsilon}
\def\i{\iota}
\def\l{\lambda}
\def\k{\kappa}
\def\6{\partial}

\def\inin{{\rm in}}
\def\lm{{\rm lm}}
\def\co{{\rm co}}
\def\lm{{\rm lm}}

\def\supp{{\rm supp}}
\def\kx{K[x]}
\def\kxx{K[[x]]}

\def\.{$\bullet$}
\def\[{[* }
%-----------------------------------------------
%          ENCODING ALGEBRAIC POWER SERIES    VVS HERWIG FEB 24  2014
%-----------------------------------------------

%\hfill  hh feb 24 2014\big\big

%\cl{\Bf  EFFECTIVE ALGEBRAIC POWER SERIES }\big\med

\cl{\Bf  ENCODING ALGEBRAIC POWER SERIES}\big\med

\cl {\smallTimes M.E.÷ ALONSO, F.J.÷ CASTRO-JIM\'ENEZ, H.÷ HAUSER}\vs 1cm

%\ignore

\footnote {}{\eightrm  M.A.÷ acknowledges support
from MEC, 2011-22435, and UCM, Grupo  910444; F.C.-J.÷ from  MTM-2010-19336 and Feder; H.H.÷ from the Austrian Science Fund FWF, projects P-21461 and P-25652.}

%-------------------------------------------------
%            ABSTRACT
%-------------------------------------------------

\cl{\vbox{\hsize 10 cm {\it Abstract}: The division algorithm for ideals of algebraic power series satisfying Hironaka's box condition is shown to be finite when expressed suitably in terms of the defining polynomial codes of the series. In particular, the codes of the reduced standard basis of the ideal can be constructed effectively.}}\vs .9cm

%----------------------------------------------- 
%            INTRODUCTION
%-----------------------------------------------			

{\Bf 1. Introduction}\med

Let $H:\A^{n+p}_K\map \A^p_K$ be a polynomial map between affine spaces over a field $K$. Assume that $H$ satisfies at $0$ the assumption of the implicit function theorem,  \med

\cl{$\6_yH(0,0)\in{\rm Gl}_p(K)$ and $H(0,0)=0$,}\med

where $y=(y_1\to y_p)$ denote coordinates on $\A^p_K$. Then there is a unique formal power series solution $h=(h_1\to h_p)$ of the system $H(x,y)=0$ at $0$, say\med

\cl{$H(x,h(x))=0$ and $h(0)=0$.}\med

Actually, the components $h_i$ are algebraic power series in the sense that each $h_i$ satisfies a univariate polynomial equation over the polynomial ring $K[x_1\to x_n]$. Conversely, any algebraic power series $h_1$ arises in this way: There is a system of polynomial equations $H(x,y)=0$ satisfying the assumption of the implicit function theorem so that the unique solution $h$ has first component $h_1$. This is known as the Artin-Mazur theorem [AM, AMR, BCR]. The characterization allows one to encode algebraic power series by a polynomial vector $H\in K[x,y]^p$ as above. The advantage of this code in comparison with taking the minimal polynomial lies in the fact that the latter determines the algebraic series only up to conjugation, so that extra information is necessary to specify the series, typically a sufficiently high truncation of the Taylor expansion. In contrast, the polynomial code $H$ determines the series $h_1$ completely and is easy to handle algebraically.\med

Phrased more abstractly, the henselization of the localization of $K[x_1\to x_n]$ at the maximal ideal $(x_1\to x_n)$ can be realized as the direct limit of finite \'etale extensions [Ar1, Na1, BCR, BrK]. Any element $h$ of the henselization, i.e., any algebraic power series, therefore belongs to such an extension -- which, by definition, can be described by a code as above.\med

It is then natural to ask to what extent operations with algebraic power series can be expressed in terms of their code; and, if this is the case for a certain operation, what will be the respective formulation of the operation in terms of the code. 
\med

In the present article we answer this question for the division of algebraic power series and for the construction of reduced standard bases of ideals. When just considered for formal power series, the division is an infinite algorithm in the infinitely many coefficients of the series. If the involved series are algebraic and satisfy Hironaka's box condition (to be defined below, see section 3), Lafon in the principal ideal case and Hironaka in general have shown that the remainder of the division is again an algebraic series [Lf, Hi1], cf.÷ also [BCR, Thm.÷ 8.2.9, p.÷ 169]. As a consequence, the reduced standard basis of the ideal is also formed by algebraic series. This fact was used for instance by Hironaka in order to construct idealistic exponents of singularities on \'etale neighborhoods and to control the behaviour of the local resolution invariant $\nu^*$ under blowup [Hi1, Hi2, chap.÷ III].\med

The beforementioned box condition is the natural extension to the case of ideals of the notion of $x_n$-regularity of a series. It postulates the existence of a specific Rees decomposition -- namely one which is generated by monomials in an appropriate coordinate system -- of the quotient module of the power series ring factored by the given ideal, cf.÷  [Re] and section 2. Algorithms to determine Rees decompositions have been proposed by Sturmfels-White [SW]. They rely on the construction of (not necessarily reduced) standard bases.\med

Starting with a system of algebraic generators of an ideal with box condition it is not at all clear how to construct, from the polynomial codes of the generators, the codes of the algebraic series defining the reduced standard basis, or, respectively, the codes of the quotients and the remainder of the division of a given algebraic series by the ideal. This question will be the subject of the article.\med

We prove that there does exist a finite algorithm which computes the codes of the reduced standard basis, respectively of the quotient and the remainder series of a division, from the codes of the algebraic input series. For the principal ideal case, i.e., the Weierstrass division, such an algorithm has been proposed and proven to work by Alonso-Mora-Raimondo [AMR]. This algorithm is already quite complicated. The general case, i.e., the division of one series by several series, is substantially more intricate and resisted for a long time.\med

In this paper we will present a complete answer to the problem, describing explicitly how to manipulate the codes of algebraic power series in order to perform the division in general. This, of course, reproves Lafon's and Hironaka's existential results on division, but it goes far beyond: It provides a quite precise manual of how to express algebraic operations with algebraic power series in terms of their codes. This is by no means trivial, and the resulting algorithm, when carried out in a concrete example, turns out to have high complexity (we give one explicit computation in the appendix). So for practical purposes the algorithm is of no big use. \med

But taken from a logical or operational point of view, the algorithm is very interesting. It is built on two simultaneous inductions, both on the number of variables, which resemble the induction which appears in the proof of the Artin approximation theorem [Ar2]. Coordinates in the affine space and generators of the ideals have to be chosen very carefully so as to make the argument work. But once this is done appropriately, the proofs develop quite naturally and are almost straightforward. In this sense, we are not only able to codify algebraic power series -- we know and understand how this codification mimics their manipulation in the division process.\med

Behind the curtain, there resides a finiteness principle which is ubiquituous in algebraic geometry and commutative algebra: The Noether normalization lemma, or, phrased differently, the finiteness of certain morphisms. In our context, this finiteness is first met in the notion of $x_n$-regularity of power series in the Weierstrass division, and then also in Hironaka's box condition and our concept of echelon (which is a Rees decomposition of a prescribed combinatorial type). It is the prerequisite for a subtle induction on the number of variables, but has the drawback that in the induction step one has to consider modules instead of ideals. This aggravates the notation, though modules are the natural context to work with. \med

The nicest part of our algorithm is what we call {\it virtual division}, a trick which has already appeared in various disguises in the literature, e.g.÷ in the work of Artin, Malgrange, Mora, Pfister-Popescu and Alonso-Mora-Raimondo: When dividing formal power series expand them with respect to one variable and write the coefficient series in the remaining variables as new unknown variables. If this is done with the necessary caution, the successive operations in the division of the formal power series can be carried out in terms of these virtual series and will then be {\it finite} processes. To make this approach work in reality, a precise understanding of the structure of the division algorithm is mandatory.\med

To resume and rephrase the above, our division algorithm for the codes of algebraic power series shows that the division is a finite process once you succeed to interpret certain packages of infinitely many data (i.e., coefficient series) as single objects which undergo a uniform transformation under division. The complexity of the algorithm shows that this encryption is by no means obvious. But it does exist and work.\med

The emphasis of the paper is theoretical -- actual computations become quickly unfeasable. 
We rather provide insight and methods of how to manipulate algebraic power series abstractly within finite algorithms. This may turn out to be useful in other situations where one aims at or needs finiteness assertions: passage to \'etale neighborhoods, noetherianity, semicontinuity of invariants of complete local rings, recursion theory for generating series, ...\med

%-------------------------------------------------

{\it Example.} Let us briefly explain the method in the special case of the construction of the code of the Weierstrass normal form of an $x_n$-regular power series $g(x)$ of order $d$. Assume for simplicity that $g$ is actually a polynomial, say $g(x)=G(x)\in K[x]$ (capital letters will be reserved throughout for polynomials). Introduce new variables $u_0\to u_{d-1}$ and define a polynomial $B\in K[x_n,u]$ as $B(x_n,u)=x_n^d +\sum_{j=0}^{d-1} u_j\cd x_n^j$. This is our candidate presentation for the Weierstrass normal form of $G$. It then suffices to determine (algebraic) series $u_0(x')\to u_{d-1}(x')\in K[[x']]=K[[x_1\to x_{n-1}]]$ such that the series $b(x)$ obtained from $B$ by substitution of $u_j$ by $u_j(x')$, say\med

\hs 3cm $b(x)=B(x_n,u(x'))=x_n^d +\sum_{j=0}^{d-1} u_j(x')\cd x_n^j$,\med

equals the Weierstrass normal form of $G$. Instead of constructing the series $u_j(x')$ directly, we shall develop a procedure to determine their code (in the sense described above, see section 6 for details). To do this, observe first that $x_n^d$ is the initial monomial of $G$ with respect to the lexicographic order $<_{lex}$ on $\N^n$ for which $(1,0\to 0)>\ldots >(0\to 0,1)$, i.e., the exponent of $x_n^d$ is the smallest element with respect to $<_{lex}$ of the support of $G$ (this uses that $u_j(0)=0$ since $b$ has order $d$ at $0$). The usual power series division of the monomial $x_n^d$ by $G$ with respect to this initial monomial then yields a formal power series remainder $r(x)=\sum_{j=0}^{d-1} u_j(x')\cd x_n^j$ such that $x_n^d-r(x)$ is the Weierstrass normal form of $G$. This division is in general an infinite process.\med

The key point now is to view $x_n^d$ alternatively as the {\it leading} monomial of the polynomial $B$ with respect to a suitable monomial order $<_\omega$ on $\N\times \N^d$, i.e., the exponent of $x_n^d$ becomes the {\it largest} element with respect to $<_\omega$ of the support of $B$. Indeed, just take for $<_\omega$ an order such that $u_j<\!\!<x_n$ for $j=0\to d-1$. Then $u_j\cd x_n^j <x_n^d$ for $j<d$ and hence $x_n^d$ will be the largest monomial of $G$ with respect to $<_\omega$. This now allows us to divide $G$ by $B$ polynomially with respect to the leading monomial $x_n^d$, say\med

\cl{ $G= Q\cd B + R$,}\med

with quotient a polynomial $Q$ in $K[x,u]$ and with remainder a polynomial $R$ in $K[x,u]$ of the form $R=\sum_{j=0}^{d-1} U_j(x',u)\cd x_n^j$ for some polynomial coefficients $U_j\in K[x',u]$. If $g$ were not a polynomial but just an algebraic series, one would have to take for $G$ the polynomial code of it, see section 13 for the precise procedure. This polynomial division is, of course, a finite process. A rather tedious computation then shows that the jacobian matrix $\6_uU$ of the vector $U=(U_0\to U_{d-1})\in K[x',u]^d$ with respect to the $u$-variables is invertible when evaluated at $0$. It thus defines, by the implicit function theorem, a unique vector $u(x')=(u_0(x')\to u_{d-1}(x'))$ of algebraic series $u_j(x')$ such that $U(x',u(x'))=0$. This just means that $U$ is a code for $u(x')$. But, by construction, $R(x,u(x'))=0$, so that $G(x)= Q(x,u(x'))\cd B(x,u(x'))$. By comparison of the initial monomials it follows that $Q(x,u(x'))$ is invertible as a power series, hence $b(x)= B(x,u(x'))$ is indeed the Weierstrass normal form of $G$ as required.\med

This example gives an idea of how the codes of reduced standard bases and of the quotients and the remainder of a division can be constructed. In practice and for the required generality the technicalities become unfortunately much more involved. \med

At the same time, there remain puzzling mysteries when the involved algebraic series are no longer $x_n$-regular (in which case the Weierstrass normal form has to be defined as the reduced  standard basis of the ideal). For instance, the polynomial $xy-z(x+y+x^2y^2)$ with initial monomial $xy$ has an algebraic series as its normal form, whereas the normal form of $xy-z(x^2+y^2+x^2y^2)$ is a transcendent series (over a ground field of characteristic zero; it is a so called Mahler series). Both facts are easy to prove by direct computation. In contrast, the normal form of $xy-z(1+y)(1+x^2y)$, a polynomial which appears in the counting of Gessel walks, is again an algebraic series, but this seems to be very intricate to prove. The algebraicity of the normal form was eventually shown by Bostan-Kauers -- a substantial part of their proof relies on heavy computer machinery [BK]. No conceptual systematic proof of the algebraicity seems to be known for this example. However, modifying slightly the input polynomial, taking now $xy-z(1+y)(1+xy^2)$, it is almost immediate to detect the algebraicity using a suitable division. \med

These examples suggest that there are hidden structural patterns which cause the phenomena to happen and which should explain the occurrence of algebraic or transcendent normal forms. Little seems to be known in this respect. For instance, the classification of the generating functions of lattice walks in the first quadrant, studied among others by Bousquet-M\'elou, Mishna and Petkov\v sek, does not seem to reveal a systematic background  [BM, BP2, Mi].\med

%-------------------------------------------------

{\it Organization of the paper.} After some preliminary recalls on the formal power series and polynomial division covering sections 2 to 5, we introduce and study in sections 6 to 8 codes of algebraic series and of the ideals generated by them. These are polynomial data which completely determine the series and ideals they encode. For later purposes the codification is carried out from the beginning for vectors of algebraic series and the modules they generate. \med

Section 9 describes how to compute the codes of standard bases of ideals and modules from a given (arbitrary) generator system (Theorem 9.1). This is straightforward, and based on Lazard's homogenization method, respectively Mora's tangent cone algorithm. Both were refined and extended by Gr\"abe and Greuel-Pfister. Our two main results (sections 10 and 11) concern the construction -- in terms of the defining codes -- of {\it reduced} standard bases of modules of algebraic power series vectors (Theorem 10.1), and of the quotients and the remainder of an algebraic power series division (Theorem 11.1).\med

The proofs of these two theorems are mutually interwoven (sections 12 to 15). First, the construction of the reduced standard basis is performed in the $x_n$-regular case (i.e., in the case where the initial module of the given module of algebraic power series vectors is generated by monomial vectors depending only on the last variable $x_n$). This is by far the most complicated step. It clearly shows how important it is to codify the series in a very systematic manner. Otherwise it would be hopeless to prove that the resulting polynomial vectors represent again codes (i.e., satisfy the assumption of the implicit function theorem). In the case of principal ideals, the proof provides the code of the Weierstrass normal form of the given series, cf.÷ [AMR].\med

The preceding construction of the codes of the reduced standard basis in the $x_n$-regular case is then used to establish the division of algebraic series on the level of codes in the $x_n$-regular case. This is not too difficult. It relies on the effectivity of the division algorithm in localizations of polynomial rings, proven by Lazard, Mora, Gr\"abe and Greuel-Pfister. \med

Once the two theorems are established in the $x_n$-regular case, the general case is carried out by induction on the number of variables. It is here that Hironaka's box condition comes into play. One key feature is its persistence under taking hyperplane sections (in a well defined sense), and this is used to know that the associated modules in $n-1$ variables satisfy again the box condition. So induction applies to prove both theorems simultaneously. \med

In the last section, we illustrate the instances and the complexity of the two algorithms in the computation of a concrete example.\med

{\it Acknowledgements.} The authors are very indebted to T.÷ Mora, W.÷ Seiler, G.-M.÷ Greuel, G.÷ Rond and O.÷ Villamayor for many valuable comments and helpful suggestions during the preparation of this article. W.÷ Seiler pointed out an inaccuracy in an earlier version of the manuscript, indicated the relation of echelons and Janet bases with his notion of $\delta$-regularity and Pommaret division, and provided several important references. Part of this work has been done during visits of the first two authors to the University of Vienna and the Erwin Schr\"odinger Institute.\med
\goodbreak

%-------------------------------------------------
%      OLD INTRODUCTION FROM JULY 27 2011
%-------------------------------------------------

\ignore

Algebraic power series are formal power series $h(x)=\sum_{\a \in \N^n} c_\a x^\a$ in several variables $x=(x_1\to x_n)$ with coefficients in a field $K$ which satisfy an algebraic relation of the form\med

\hs 2cm $p_d(x)h(x)^d+p_{d-1}(x)h(x)^{d-1}+\ldots +p_1(x)h(x)+p_0(x)=0$,\med

where the coefficients $p_i(x)$ are polynomials. Typical examples are rational  functions as $x\cd(1+x)\inv$, roots of polynomials as $\sqrt {1+x^2y}$, inverses $f\inv$  of polynomial mappings $f:K^n\map K^n$ satisfying at a point $p$ the assumption of the Inverse Function Theorem as $f(x,y)=(x+x^3,y-xy^2)$ at $0$, or solutions $y(x)$ of polynomial equations $f(x,y)=0$ satisfying at a point $p$ the assumption of the implicit function theorem w.r.t.÷ the variables $y$ as $f(x,y)=y+xy+x^3y^2$ at $0$.\med

The ring of algebraic series in $n$ variables is thus the algebraic closure of the polynomial ring $K[x_1\to x_n]$ inside the formal power series ring $K[[x_1\to x_n]]$. It can equivalently be interpreted as the henselization of the polynomial ring at $0$ [Ar1, Ra, BrK, p.÷ 333]. \med

Note that the minimal polynomial of an algebraic series $h$ determines $h$ only up to conjugacy: there may be other power series solutions to the equation, the conjugates of $h$, and $h$ can be distinguished from these for instance by a sufficiently high truncation of its power series expansion. The simplest example thereof is the equation $y^2 - 2y + x=0$ with algebraic solutions $h_{\pm}= 1\pm\sqrt{1-x}$.\med 

There is an alternative way to encode algebraic series, using the Artin-Mazur theorem [AM, p.÷ 88, AMR, appendix, BCR, Thm.÷ 8.4.4, p.÷ 173]: the algebraic series then appears as the first component of the solution vector of a system of polynomial equations in two sets of variables satisfying the assumption of the implicit function theorem at $0$ as in the last example above. This type of code is more convenient to handle when performing algebraic operations with the series. Our interest here will be to use this code specifically for the power series division of algebraic series. \med

The classical Weierstrass Division Theorem for formal or convergent power series was extended by Grauert, Hironaka and Galligo to the case where a series is divided by several series, say an ideal [Gra, Hi1, Ga]. The remainder series is made unique by imposing conditions on the allowed monomials, which are required not to be divisible by any of the initial monomials of the elements of the ideal. Lafon in the principal ideal case and Hironaka in general showed that the remainder series is algebraic if the input series were algebraic, provided that the initial ideal satisfies a certain combinatorial condition -- the so called box condition -- which extends the notion of $x_n$-regularity of Weierstrass [Lf, Hi1].  
Without any condition on the ideal the algebraicity may fail [Hi1, p.÷ 75]. \med

The proofs of Lafon and Hironaka are existential -- they do not indicate a method how to find the minimal polynomial or the code of the remainder series. The formal power series division theorem allows us to compute the division of formal power series up to arbitrarily high degree. If the remainder is known to be algebraic, its expansion up to large degree does not, however, suffice to determine its code.  Actually, it is not at all clear how to compute effectively the remainder in the case that all involved series are algebraic. \med

In the Weierstrass case, i.e., the division of a series by {\it one} $x_n$-regular series, this effectivity problem was settled affirmatively in 1992 by Alonso, Mora and Raimondo [AMR]. Their answer is rather tricky and involved, and gives, among other things, an algorithm to compute the code of the Weierstrass form of the series. The method does not extend directly to the case of the division by an ideal with {\it several} generators.\med

The objective of the present paper is to show that it is indeed possible to perform effectively   in terms of the respective codes the power series division of an algebraic series by an ideal of algebraic series satisfying the box condition: there is a {\it finite algorithm} to compute the code of the output series from the input codes. The algorithm is highly involved, using a multiple induction on the generators of the ideal and the number of variables, and is thus only theoretically feasable. The point is that it exhibits precise rules of how to manipulate ideals of algebraic series in an explicit and constructive manner. \med

To do so, we first have to develop a suitable framework and language to work with codes of algebraic power series. This concerns in particular the construction of codes for ideals and modules of algebraic series. Then, simple procedures as Mora's tangent cone algorithm and Lazard's homogenization method have to be expressed on the level of codes. This is a prerequisite to compute by a finite algorithm a standard basis of an ideal of $K[[x]]$ generated by polynomials (these polynomials will be in the application the codes of algebraic series). \med

The hardest part in rendering the algebraic series division effective is the construction of a {\it reduced} standard basis of an ideal of algebraic series. The clue for this is the notion of a {\it virtual} reduced standard basis. This is a system of polynomials with unknown coefficients which anticipates the actual form of the reduced standard basis. The unknowns can then be determined by applying the polynomial division algorithm of Gr\"obner bases theory. The resulting coefficients are not given through their power series expansions but -- again -- through polynomial codes. In particular, they are shown to be algebraic series. Both sets of polynomials  together, the virtual reduced standard basis and these codes, then define the code for the {\it actual} reduced standard basis.
\med

Once the code for the reduced standard basis of the ideal is found, the division of an algebraic series by this basis goes by induction on the number of variables. The procedure uses only the codes and is rather straightforward.\med

We emphasize that the various constructions of this paper are not invented for practical purposes. Actually, it seems that a low complexity algorithm for the division of algebraic series cannot even exist. In contrast, the systematic encoding of algebraic power series and the explicit manipulation of the codes provide techniques which may allow to work effectively with algebraic series also in other situations.\med

\recognize

\big \goodbreak

%----------------------------------------------- 
%            PRELIMINARIES
%-----------------------------------------------			

%\cl {\Bf PRELIMINARIES}\big

%-------------------------------------------------
%          MONOMIAL MODULES
%-------------------------------------------------

{\Bf 2. Monomial modules}\med

The letters $n$, $p$, $r$, $s$ are reserved for fixed integers in $\N$. The
letters $i$, $k$ and $\ell$  will generally vary in the ranges  $1\leq i\leq p$,  $1\leq
k\leq r$ and $1\leq \ell\leq s$.\med

We denote by $K[x_1\to x_n]=\kx$ and $K[[x_1\to x_n]]=\kxx$ the polynomial,
respectively formal power series ring in $n$ variables $x=(x_1\to x_n)$ over a
field $K$.  Elements of $\kx^s$ and $\kxx^s$ will be
called {\it polynomial} respectively {\it formal power series vectors}. Capital
letters will be reserved for polynomials,  lower case letters for power series. 
We set $x'=(x_1\to x_{n-1})$ and denote by $y=(y_1\to y_p)$ additional
variables.\med

Vectors $g\in \kxx^s$ will be expanded into $g= \sum_{\a \ell} c_{\a \ell}x^\a
e_\ell$ with $c_{\a \ell}\in K$ and $e_\ell=(0\to 0,1,0\to 0)$ the
canonical $K$-basis of $K^s$.  The vectors $x^\a e_\ell$ are called {\it monomial vectors}. Note that all their entries but one are zero: a vector all whose entries are monomials will not be considered here as a monomial vector. The {\it support} of $g$ is the set   $\supp(g)=\{(\a,\ell)\in \N^n\times \{1\to s\},\, c_{\a \ell}\neq  0\}$. We sometimes abbreviate pairs $(\a,\ell)$ by $\a\ell$.   \med

Brackets $\<g_1\to g_r\>$ denote submodules of $K[[x]]^s$ generated by
power series vectors $g_1\to g_r\in \kxx^s$. We abbreviate this by
$\<g_k\>$ if the range of $k$ is clear from the context. \med

A {\it monomial submodule} of $\kxx^s$ is a submodule $M$ of $\kxx^s$ generated by monomial vectors. It is a cartesian product $M=\Pi_{\ell=1}^s M_\ell$ of monomial ideals
$M_\ell$ in $\kxx$. The elements of $M$ are the power series vectors with support
in $\Sigma=\{(\a,\ell)\in \N^n\times \{1\to s\},\, x^\a e_\ell\in M\}$. The {\it
canonical direct monomial complement} of a monomial submodule $M$ of $\kxx^s$ is the subvectorspace $\co(M)$ of $\kxx^s$ of power series vectors with support in $\Sigma'=(\N^n\times \{1\to s\})\sm \Sigma$. This provides the direct sum decomposition of $K$-vectorspaces $M\oplus \co(M)=\kxx^s$.\med

We say that a monomial submodule $M$ of $K[[x]]^s$ is $x_n$-{\it regular} if it is
generated by monomial vectors in $K[[x_n]]^s$, say $M=\<M\cap K[[x_n]]^s\>$. We shall then always assume for simplicity -- applying
if necessary a permutation of the components of $\kxx^s$ -- that $M$ is
generated by vectors of the form $x_n^{d_k}\cd e_k$ with $d_k\geq 0$ and $1\leq
k\leq r$ for some $r\leq s$. In this case the complement $\co(M)$ is a cartesian product\med

\hs 3cm $ \co(M)=\prod_{k=1}^{r} (\oplus_{j=0}^{d_k-1} K[[x']]\cd x_n^j) \times\kxx^{s-r}$\med

of a finitely generated free $K[[x']]$-module with a finitely generated free $K[[x]]$-module. We say that $M$ satisfies {\it Hironaka's box condition} if $\co(M)$ can be written as a cartesian product of direct sums of finite free monomial $K[[x_1\to x_j]]$-modules\med

\hs 3cm$\co(M) =\prod_{\ell=1}^s\, \oplus_{j=0}^n\, \oplus_{\gamma\in \Gamma_{\ell j}} K[[x_1\to x_j]]\cd x^\gamma$\med\goodbreak

with finite sets $\Gamma_{\ell j}\subset \N^n$. Being $x_n$-regular is a special case of the box condition. For cyclic submodules of $\kxx^s$, both notions coincide. They obviously depend on the numbering of the variables $x_1\to x_n$. Notice that for $s=1$ and $0\neq M\subsetneq K[[x]]$ a non trivial ideal, the indices of the boxes $F_j$ run from $1$ to $n-1$. Also notice that the box condition for a monomial submodule $M\subset K[[x]]^s$ is equivalent to the box condition for each of the factors of $M$ (which are monomial ideals in $K[[x]]$).\med

W.÷ Seiler informed us that in the case of ideals the box condition is equivalent to his notion of $\delta$-regular coordinates [Se3]. We say that a monomial submodule $M$ of $K[[x]]^s$ is an {\it echelon} if it can be written as\med

\hs 3cm$M=\prod_{\ell=1}^s\, \oplus_{j=0}^n\, \oplus_{\delta\in \Delta_{\ell j}}
K[[x_1\to x_j]]\cd x^\delta$\med

with finite sets $\Delta_{\ell j}\subset \N^n$. This can be rewritten as \med

%\hs 3cm$M=\prod_{\ell=1}^s\, \oplus_{\delta\in \Delta_\ell}\ K[[x_1\to x_{n_\delta}]]\cd x^\delta$\med

\hs 3.35cm $M = \oplus_{\ell=1}^s\, \,\oplus_{\delta\in \Delta_\ell}K[[x_1\to x_{n_\delta}]]\cd x^\delta\cd e_\ell$\med

where $\Delta_\ell =\bigcup_j \Delta_{\ell j}$ and where, for each $\delta$, the index $n_\delta$ takes a value between $0$ and $n$. This notion is a special case of a Rees decomposition of $M$ [Re]. We call the collection of monomial vectors $x^\delta\cd e_\ell$ with $\delta\in \Delta_\ell$ and $1\leq \ell \leq s$ a {\it Janet basis} of the echelon $M$ with {\it scopes} $n_\delta$ (also known as {\it levels} or {\it classes}). Our definition differs slightly from Janet's original definition in the sense that we only allow nested groups of variables in the coefficients [Ja1, Ja2], cf.÷ also [Ri]. We refer to the related notions of Pommaret bases and involutive bases [GB, Se1, Se2], and the more general concepts of Rees and Stanley decompositions of rings [Re, SW, Am, BG].  \med
\goodbreak

%-------------------------------------------------

For the following result, see also Janet [Ja1, Ja2] and Seiler [Se2]. \med

{\bf Theorem 2.1.} {\it Monomial submodules of $K[[x]]^s$ satisfying Hironaka's box conditon are echelons.}\med

%-------------------------------------------------

{\it Proof.} Let $M$ be such a module, and let $M_n=\<M\cap K[[x_n]]^s\>$ be the submodule of $K[[x]]^s$ generated by the $x_n$-pure monomial vectors of $M$. By definition, $M_n$ is $x_n$-regular. Let $x_n^{d_k}\cd e_k$ with $1\leq k \leq r$ be a minimal generator system of $M_n$ (after possibly permuting the components of $K[[x]]^s$). Then $M_n=\oplus_{k=1}^r K[[x]]\cd x_n^{d_k}\cd e_k$, which shows that the monomial vectors $x_n^{d_k}\cd e_k$ form a Janet basis of $M_n$ with scopes $n_k=n$. The direct sum decomposition\med

\hs 1cm $K[[x]]^s=M_n \oplus  (\oplus_{m=1}^r\,\oplus_{j=0}^{d_m-1} K[[x']]\cdot x_n^j\cd e_m) \,\, \oplus\,\,\oplus_{m=r+1}^s\, K[[x]]\cd e_m$\med

yields a decomposition $M=M_n\oplus M'$ where $M'$ is now a $K[[x']]$-submodule of the  finitely generated free $K[[x']]$-module $\oplus_{m=1}^r\,\oplus_{j=0}^{d_m-1} K[[x']]\cdot x_n^j\cd e_m$. We use here that, because of the box condition, $M$ has zero intersection with $\oplus_{m=r+1}^s\, K[[x]]\cd e_m$. \med

It is checked that the box condition persists under the above decomposition, i.e., that $M'$ satisfies it again. By induction on the number of variables, $M'$ is an echelon. Its Janet basis has scopes $\leq n-1$. From $M=M_n\oplus M'$ now follows that also $M$ is an echelon.\med

%-------------------------------------------------

{\it Example.} The assertion of the theorem does not hold for arbitrary modules as was pointed out by W.÷ Seiler. Take the ideal $I$ of $K[x,y,z]$ generated by the three monomials $xy$, $xz$ and $yz$. It is easy to see that it does not satisfy the box condition. And it is not an echelon, since, for instance, among the monomials of $I$ which are not multiples of $xy$ one has monomials $x^dz$ and $y^dz$ of arbitrary degree $d$ in $x$ and $y$. As the situation is symmetric with respect to any permutation of the variables, $I$ does not admit the required decomposition of an echelon. \med
\big

\ignore
%-------------------------------------------------
%       EXTRACT FROM ECHELON PREPRINT AUGUST 11 2004
%-------------------------------------------------

The next two lemmata are slightly more involved. All assertions hold correspondingly for
ideals and echelons in $\zeta+\N^r$ with $\zeta\in \Z^n$ and $r\leq n$.\med

{\bf Lemma 2}.  (a) {\it The intersection of an ideal, resp.÷ echelon $E$ in
$\N^n$ with an azulejo $\d+\N^k$ in $\N^n$ is an ideal, resp.÷ echelon  in
$\N^n\cap(\d+\N^k)$.}

(b) {\it Ideals are echelons. More explicitly, for any ideal $E$ in $\N^n$, there exists a mosaique $\P=\{A_\a\}_{\a\in V}$ in $\N^n$ with support
$E$, $E=\dot\cup_{\a\in V}\, A_\a$.}

(c) {\it If an echelon $E$ satisfies Hironaka's box condition in $\N^n$, its intersection
$E'=E\cap (\d+\N^k)$ with any azulejo $\d+\N^k$ in $\N^n$ satisfies the box
condition in $\N^n\cap(\d+\N^k)$.}\med

%-------------------------------------------------------

{\it Proof.} For (a), it suffices to recall that by Lemma 1(a) the intersection of azulejos is
again an azulejo. For (b), write $E$ as a finite (not necessarily disjoint)
union $E=\cup_{\a\in V}\, \a+\N^n$ with $V\subset E$ finite. Let $d$ be the
maximum of the last components $\a_n$ of elements of $V$. Then $E$ is a
finite disjoint union of the intersections $E_i=E\cap (\N^{n-1}\times\{i\})$
for $0\leq i\leq d-1$ and $E(d)=E\cap ((0\to 0,d)+\N^n)$. The sets $E_i$ are
by (a) ideals in $\N^{n-1}\times\{i\}$, which, by induction on $n$, are
echelons. As for $E(d)$, observe that all $\a\in V\cap E(d)$ have last
component $\a_n=d$. Hence $E(d)$ is a cartesian product $\tilde E \times
(d+\N)$ for some ideal $\tilde E$ in $\N^{n-1}$. Using again induction on
$n$, $\tilde E$ admits a mosaique of azulejos, hence also $E(d)$. This proves
(b).\med

For (c), let $F=\N^n\setminus E$ be the complement of $E$. Then \med

\hs 1cm $F'=F\cap (\d+\N^k) = (\N^n\setminus E) \cap 
(\d+\N^k) =(\d+\N^k)\setminus E= (\d+\N^k)\setminus E'$ \med

is the complement of $E'$ in $\d+\N^k$. The assumption on $E$ implies that 
$F$ admits a partition $F=\dot\cup_{\gamma\in Z}\, (\gamma+\N^{n_\gamma})$ in
azulejos. By assertion (a), $E'$ is an echelon. Hence it is sufficient to
construct a partition of $F'$ in azulejos. But by Lemma 1 (a), the intersections
$(\gamma+\N^{n_\gamma})\cap (\d+\N^k)$ are azulejos, so that (c) follows from
\med

\hs 1cm $F'= (\d+\N^k)\cap F=(\d+\N^k)\cap (\dot\cup_{\gamma\in Z}\,
(\gamma+\N^{n_\gamma}))=\dot\cup_{\gamma\in Z}  ((\d+\N^k)\cap
(\gamma+\N^{n_\gamma}))$. \med

%-------------------------------------------------------

{\bf Lemma 3}.  (a) {\it If a non-empty ideal $E$ has mosaique $E=\dot\cup_{\a\in
V}\, (\a+ \N^{n_\a})$ and satisfies Hironaka's box condition in $\N^n$, there is
precisely one $\a\in V$ which is of form $\a=(0\to 0,d)$ for some $d\geq 0$
and so that $n_\a=n$.} 

(b) {\it If an echelon $E$ has mosaique
$E=\dot\cup_{\a\in V}\, (\a+ \N^{n_\a})$ and satisfies Hironaka's box
condition in $\N^n$ with complement $F=\N^n\sm E=\dot\cup_{\gamma\in Z}\,
(\gamma+ \N^{n_\gamma})$, there is precisely one $\delta\in V\cup Z$ of form
$\delta=(0\to 0,d)$ with scope $n_\delta=n$.} \med

%-------------------------------------------------------

{\it Proof.} As for assertion (a), there can be at most one $\a\in V$ of scope $n$, for the
azulejos of the mosaique of $E$ are disjoint. Assume there is none. Then $F$
would have elements $\gamma\in\N^n$ with arbitrarily large last component
$\gamma_n$. Hence the mosaique of $F$ would have one azulejo of scope $n$.
But then, since $E$ is an ideal and not empty, $E\cap F\neq \emptyset$,
contradiction. So there is precisely one $\a\in V$ of scope $n$. It is then
clear that this $\a$ must be of form $\a=(0\to 0,d)$.\med

Assertion (b) goes similarly. If
all scopes $n_\a$ of $\a\in E$ are $<n$, at least one $\delta\in Z$ must have
scope $n$, else $E\cup F\subsetneq \N^n$.  As $F=\dot\cup_{\gamma\in Z}\,
(\gamma+ \N^{n_\gamma})$ is a disjoint union, there is precisely one such
$\delta$. Using again $E\cup F= \N^n$ we get $\delta=(0\to 0,d)$ for some
$d\in \N$.\med

%-------------------------------------------------------

{\bf Lemma 3}.  (a) {\it If a non-empty ideal $E$ has mosaique $E=\dot\cup_{\a\in
V}\, (\a+ \N^{n_\a})$ and satisfies Hironaka's box condition in $\N^n$, there is
precisely one $\a\in V$ which is of form $\a=(0\to 0,d)$ for some $d\geq 0$
and so that $n_\a=n$.} 

(b) {\it If an echelon $E$ has mosaique
$E=\dot\cup_{\a\in V}\, (\a+ \N^{n_\a})$ and satisfies Hironaka's box
condition in $\N^n$ with complement $F=\N^n\sm E=\dot\cup_{\gamma\in Z}\,
(\gamma+ \N^{n_\gamma})$, there is precisely one $\delta\in V\cup Z$ of form
$\delta=(0\to 0,d)$ with scope $n_\delta=n$.} \med

%-------------------------------------------------------

{\it Proof.} As for assertion (a), there can be at most one $\a\in V$ of scope $n$, for the
azulejos of the mosaique of $E$ are disjoint. Assume there is none. Then $F$
would have elements $\gamma\in\N^n$ with arbitrarily large last component
$\gamma_n$. Hence the mosaique of $F$ would have one azulejo of scope $n$.
But then, since $E$ is an ideal and not empty, $E\cap F\neq \emptyset$,
contradiction. So there is precisely one $\a\in V$ of scope $n$. It is then
clear that this $\a$ must be of form $\a=(0\to 0,d)$.\med

Assertion (b) goes similarly. If all scopes $n_\a$ of $\a\in E$ are $<n$, at least one $\delta\in Z$ must have scope $n$, else $E\cup F\subsetneq \N^n$.  As $F=\dot\cup_{\gamma\in Z}\, (\gamma+ \N^{n_\gamma})$ is a disjoint union, there is precisely one such
$\delta$. Using again $E\cup F= \N^n$ we get $\delta=(0\to 0,d)$ for some
$d\in \N$.\med

\big\goodbreak

%-------------------------------------------------
%   END EXTRACT FROM ECHELONS
%-------------------------------------------------

\recognize

%----------------------------------------------- 
%           MONOMIAL ORDERS AND INITIAL MODULES
%-----------------------------------------------			

{\Bf 3. Monomial orders and initial modules}\med

Division theorems are based on ordering the summands $c_{\a \ell}x^\a e_\ell$ of the expansion of a power series vector $g= \sum_{\a \ell} c_{\a \ell}x^\a e_\ell$ according to the indices $(\a,\ell)\in \N^n\times \{1\to s\}$ with non-zero coefficients $c_{\a \ell}$: A {\it monomial order} on $\N^n\times \{1\to s\}$ is a total order $<_\eta$ on $\N^n\times \{1\to s\}$ which is compatible with the semi-group structure of $\N^n$, having $0$ as its smallest element, and which is noetherian. This means that  if $(\a,\ell)<_\eta (\b,m)$ then  $(\a+\gamma,\ell)<_\eta (\b+\gamma,m)$ for any $\gamma\in\N^n$, and, secondly, that any decreasing sequence becomes stationary. The order is {\it degree
compatible} if $\abs\a<\abs\b$ implies $(\a,\ell)<_\eta(\b,m)$, where $\abs\a$ denotes the sum of the components of $\a$. An {\it extension} 
of $<_\eta$ is a monomial order $<_\eps$ on $\N^{n+p}\times \{1\to s\}$ whose
restrictions to $\N^n\times\{\delta\}\times \{1\to s\}$ coincide  for all $\delta\in\N^p$
with the order induced by $<_\eta$ on $\N^n\times\{\delta\}\times \{1\to s\}$. We will always identify monomial orders on $\N^n\times \{1\to s\}$ with the
induded ordering of the monomial vectors in $\kxx^s$.\med

The {\it initial monomial vector} $\inin(g)$ of  $g=\sum c_{\a \ell}x^\a 
e_\ell\in \kxx^s$ with respect to $<_\eta$ is the vector $x^\a \cd e_\ell$ of
the expansion of $g$ for which $(\a,\ell)$ is {\it minimal} with respect to $<_\eta$. We shall assume that $x^\a \cd e_\ell$ has coefficient $1$ in the expansion of $g$. We
then write   $g=x^\a \cd e_\ell-\ol g$ and call $\ol g$ the {\it tail} of $g$.
\med  

For a submodule $I$ of $\kxx^s$, the {\it initial module} of $I$ with respect  to
$<_\eta$ is the monomial submodule $\inin(I)$ of $\kxx^s$ generated by all initial
monomial vectors of elements of $I$. This is a monomial submodule which depends on the choice of $<_\eta$. We denote by $\co(I)$ the canonical direct monomial complement of $\inin(I)$ in $\kxx^s$. Elements $g_1\to g_r$ of $\kxx^s$ form a {\it standard basis} w.r.t.÷ $<_\eta$ if their initial monomial vectors generate the initial module $\inin(I)$ of the module $I$ generated by $g_1\to g_r$. They are a {\it reduced standard basis} if the tails $\ol g_k$ belong to  $\co(I)$.  We do not require that a reduced standard basis is minimal.\med

We say that a submodule $I$ of $K[[x]]^s$ is $x_n$-{\it regular}, respectively satisfies {\it Hironaka's box condition}, or is an {\it echelon} with respect to the monomial order $<_\eta$ on $\N^n\times\{1\to s\}$, if its initial module $\inin(I)$ is $x_n$-regular, respectively satisfies the box condition, or is an echelon. A {\it Janet basis} of a submodule $I$ of $K[[x]]^s$ which is an echelon w.r.t.÷ $<_\eta$ is a generator system $g_1\to g_r$ of $I$ whose initial monomial vectors $\inin(g_k)$ form a Janet basis of $\inin(I)$.
\med
For a polynomial vector $G\in K[x]^s$, define the {\it leading monomial vector} $\lm(G)$ as the monomial vector $x^\a e_\ell$ of the expansion of $G$ which is {\it maximal} with respect to the chosen monomial order. Similarly as for initial modules, one obtains now the leading module $\lm(I)$ of a submodule $I$ of $K[x]^s$.\med
\big 

%-----------------------------------------------
%     DIVISION THEOREM FOR FORMAL POWER SERIES
%-----------------------------------------------

{\Bf 4. Division of formal power series and polynomials}\med

We recall the division theorem for modules of formal power series of
Grauert, Hironaka and Galligo [AHV, Gra, Hi1, Ga, HM]. For extensions of this result to more general settings see [Am, BG, GB, Se2]. \med\goodbreak

%-----------------------------------------------

{\bf Theorem 4.1.} {\it Let $I$ be a submodule of $\kxx^s$ with initial module
$\inin(I)$ with respect to a monomial order $<_\eta$ on $\N^n\times \{1\to s\}$.
Let $\co(I)$ be the canonical direct monomial complement of $\inin(I)$ in
$\kxx^s$. Then $I\oplus \co(I)=\kxx^s$.}\med
\goodbreak
%-----------------------------------------------

{\it Sketch of proof.} The sum $I\oplus \co(I)$ is direct by definition of $\co(I)$. To see that it equals $K[[x]]^s$, choose a standard basis $g_1\to g_r$ of $I$. It suffices to show that the linear map $u: K[[x]]^r\times \co(I)\map K[[x]]^s,\, (a_1\to a_r,b)\map \sum a_kg_k +b$ is surjective. By definition of standard bases, the map $v: K[[x]]^r\times \co(I)\map K[[x]]^s,\, (a_1\to a_r,b)\map \sum a_k\cd\inin(g_k) +b$ is surjective. Writing $u=v+w$, the assertion follows by restricting $v$ to a direct complement $L$ of its kernel and by showing  that $u_{\vert L}= v_{\vert L} + w_{\vert L}$ is an isomorphism with inverse the geometric series $(v_{\vert L})\inv\sum_{j=0}^\infty ((v_{\vert L})\inv w_{\vert L})^j$. This series then induces the required linear map $K[[x]]^s\map L$ inverse to $u_{\vert L}$, see [HM, Thm.÷ 5.1] for details.\med

The division theorem can be formulated more explicitly as follows: If $g_1\to g_r$ generate $I$, each vector $f\in\kxx^s$ has a decomposition $f=\sum_k a_kg_k+ h$ with unique $h\in \co(I)$. The power series expansions of the quotients $a_k$ and the remainder $h$ can be computed up to any given degree by a finite algorithm (take the expansion of the geometric series above up to the respective degree). The requirement that $h$ belongs to $\co(I)$ makes the remainder independent of the choice of $g_1\to g_r$ (but it depends on the monomial order $<_\eta$). If $g_1\to g_r$ form a standard basis, the quotients $a_k$ can be made unique by imposing suitable support conditions on them [Ga]. A reduced standard basis of $I$ is given as $x^{\a}\cd e_{\ell}-h_{\a\ell}$ with $(\a,\ell)$ varying in some finite subset $V\subset \N^n\times\{1\to s\}$, where the vectors $x^{\a}\cd e_{\ell}$ are generators of $\inin(I)$ and the vectors $h_{\a\ell}$ denote the
remainder of the division of $x^{\a}\cd e_{\ell}$ by $I$.\med 

%-----------------------------------------------

For modules which are echelons one can formulate a more precise statement:\med

{\bf Theorem 4.2.} {\it Let $I$ be a submodule of $\kxx^s$ with initial module
$\inin(I)$ w.r.t.÷ a monomial order $<_\eta$ on $\N^n\times \{1\to s\}$. Assume that $I$ is an echelon, and let $x^\a\cd e_\ell$ be a Janet basis of $ \inin(I)$ with scopes $n_{\a\ell}$, $(\a,\ell)$ varying in some finite set $V\subset\N^n\times\{1\to s\}$. Choose any elements $g_{\a\ell}$ of $I$ with initial monomial vectors $x^\a\cd e_\ell$. Then \med

\hs 2cm $ I\oplus \co(I)=\oplus_{\a\ell\in V}\ K[[x_1\to x_{n_{\a\ell}}]]\cd g_{\a\ell}\oplus \co(I)=\kxx^s$.}\big 

%-----------------------------------------------

{\it Proof.} First notice that $\inin(I)=\oplus_{\a\ell\in V}\ K[[x_1\to x_{n_{\a\ell}}]]\cd x^\a\cd e_\ell$ by definition of echelons. This allows us to modify the map $u$ from the proof of the division theorem by restricting it to the $K$-subspace \med

\hs 3cm $\oplus_{\a\ell\in V}\ K[[x_1\to x_{n_{\a\ell}}]]\cd x^\a\cd e_\ell \times \co(I)$.\med

The map $v$ is then by construction an isomorphism, and the same reasoning as before shows that this holds also for $u$. This proves the claim. \med

%-------------------------------------------------

In the polynomial case, the division admits an analogous formulation. The same proof as above applies, because the evaluation of the geometric series $(v_{\vert L})\inv\sum_{j=0}^\infty ((v_{\vert L})\inv w_{\vert L})^j$ on a polynomial vector $(a_1\to a_r,b)\in K[x]^r\times \co(I)$ terminates at sufficiently large $j$. \med

{\bf Theorem 4.3.} {\it Let $I$ be a submodule of $\kx^s$ with leading module
$\lm(I)$ with respect to a monomial order $<_\eta$ on $\N^n\times \{1\to s\}$.
Let $\co(I)$ be the canonical direct monomial complement of $\lm(I)$ in
$\kx^s$. Then $I\oplus \co(I)=\kx^s$.}\med

%-------------------------------------------------

Again, there is a more precise version in case the leading module $\lm(I)$ is an echelon. \med

{\bf Theorem 4.4.} {\it Let $I$ be a submodule of $K[x]^s$ with leading monomial module $\lm(I)$ with respect to a monomial order $<_\eta$ on $\N^n\times \{1\to s\}$. Assume that $\lm(I)$ is an echelon. Let $G_k$ be a polynomial Janet basis of $I$ with leading monomial vectors $\lm(G_k)$ of scope $n_k$. Then any $F\in K[x]^s$ admits a unique division \med

\hs 4cm $F=\sum_k A_k G_k + C$ \med

with $A_k\in K[x_1\to x_{n_k}]$ and $C\in \co(I)$. The decomposition can be obtained from the polynomial vectors $F$ and $G_k$ by a finite algorithm.}\med

\big\goodbreak

%-----------------------------------------------
%     ALGEBRAIC POWER SERIES
%-----------------------------------------------

{\Bf 5. Algebraic power series }\med

Algebraic power series are formal power series $h(x)=\sum_{\a \in \N^n} c_\a x^\a$ in several variables $x=(x_1\to x_n)$ with coefficients in a field $K$ which satisfy an algebraic relation of the form\med

\hs 2cm $P(x,h(x))=p_dh^d+p_{d-1}h^{d-1}+\ldots +p_1h+p_0=0$,\med

where the coefficients $p_i=p_i(x)$ are polynomials. We refer to [Ar2, BCR, BrK, KPR, Lf,  Na1, Na2, Ra, Ru, Wi] for the respective background. An algebraic power series vector is a vector in $\kxx^s$ whose components are algebraic series. \med

Typical algebraic series are rational functions as $x\cd(1+x)\inv$, roots of polynomials as $\sqrt {1+x^2y}$, inverses $f\inv$  of polynomial mappings $f:K^n\map K^n$ satisfying at a point $p$ the assumption of the Inverse Function Theorem as $f(x,y)=(x+x^3,y-xy^2)$ at $0$, or solutions $y(x)$ of polynomial equations $f(x,y)=0$ satisfying at a point $p$ the assumption of the implicit function theorem with respect to the variables $y$ as $f(x,y)=y+xy+x^3y^2$ at $0$.\med

The ring of algebraic series in $n$ variables is thus the algebraic closure of the polynomial ring $K[x_1\to x_n]$ inside the formal power series ring $K[[x_1\to x_n]]$. It can equivalently be interpreted as the henselization of the polynomial ring at $0$. \med

Note that the minimal polynomial of an algebraic series $h$ determines $h$ only up to conjugacy: there may be other power series solutions to the equation, the conjugates of $h$, and $h$ can be distinguished from these for instance by a sufficiently high truncation of its Taylor expansion. The simplest example thereof is the equation $y^2 - 2y + x=0$ with algebraic solutions $h_{\pm}= 1\pm\sqrt{1-x}$.\med 

%-------------------------------------------------

Lafon proved in 1965 that the Weierstrass division preserves the algebraicity of the involved series [Lf], see also [BCR]. This was reproven in 2000 by Bousquet-M\'elou and Petkov\v sek working with the recursions defining the coefficients of the series [BP1]. The result of Lafon was extended by Hironaka in 1977 to the division by ideals with several generators satisfying the box condition [Hi1]. We formulate here the division directly for modules.\med

%-----------------------------------------------

{\bf Theorem 5.1.} {\it Let $I$ be a submodule of $\kxx^s$ generated by algebraic
power series vectors. Assume that $I$ satisfies Hironaka's box condition with
respect to a monomial order $<_\eta$ on $\N^n\times \{1\to s\}$. For
any algebraic power series vector $f\in\kxx^s$ the remainder $c$ of the formal
power series division of $f$ by $I$ with respect to $<_\eta$ is an algebraic power series vector.}\med

%-----------------------------------------------

The theorem implies in particular that any submodule of $K[[x]]^s$ with box condition which is generated by algebraic power series vectors admits a reduced standard basis consisting of algebraic power series vectors. Without box condition the remainder of the division need not be algebraic. In  [Hi1, p.÷ 75], Hironaka cites the following example of Gabber and Kashiwara, which was rediscovered by Bousquet-M\'elou and Petkov\u sek in combinatorics when counting lattice paths [BP1, BP2]. \med

{\bf Example 5.2.} Divide $xy$ by $g=(x-y^2)(y-x^2)=xy-x^3-y^3+x^2y^2$ as formal power series with respect to the initial monomial $xy$. The remainder of the division lies in $\co(xy)=K[[x]]+K[[y]]$ and equals the lacunary series $b=\sum_{k\geq 0} x^{3\cd 2^k}+ \sum_{k\geq 0} y^{3\cd 2^k}$ which is transcendent . Alternatively, we may write $xy = a\cd g + r(x)+s(y)$ with series $a\in K[[x,y]]$, $r\in K[[x]]$, $s\in K[[y]]$. The symmetry between $x$ and $y$ in this expression yields $r(x)=s(x)$. Substituting $y$ by $x^2$ produces $x^3 = a\cd 0+ r(x)+r(x^2)$ which also gives the expansion of $r$. \med 

\big\goodbreak

%-------------------------------------------------
%          RESULTS
%-------------------------------------------------

%\cl {\Bf RESULTS}\big

%----------------------------------------------- 
%            CODES OF ALGEBRAIC SERIES
%-----------------------------------------------			

{\Bf 6. Codes of algebraic power series}\med

In this section we introduce the necessary terminology for working effectively with algebraic power series. The variables $x=(x_1\to x_n)$ and $y=(y_1\to y_p)$ are fixed throughout.\med

A {\it mother code} (over $x$ and $y$) is a polynomial row vector 
$H=(H_1\to H_p)\in K[x,y]^p$ with $H(0,0)=0$ whose Jacobian  matrix $D_yH$ with
respect to $y$ is invertible at $0$,\med

\hs 3cm $D_yH(0,0)\in {\rm Gl}_p(K)$.\med

The invertibility of $D_yH(0,0)$ can be rephrased by saying that for any
degree compatible monomial order on $\N^p$ the initial ideal of the ideal
$\<H_1(0,y)\to H_p(0,y)\>$ of $K[[y]]$ is generated by $y_1\to y_p$. There then exists
a linear coordinate change in the $y_i$'s so that the initial monomials
$\inin(H_i(0,y))$ of $H_i(0,y)$ equal $y_i$. For any degree compatible monomial order on $\N^{n+p}$ so that $y_i< x_j$ for all $i$ and $j$ it then follows that the initial monomials $\inin(H_i)$ of $H_i$ equal $y_i$. Instead of changing the $y_i$'s one
could also change $H$ by multiplying it from the right with a suitable matrix
in  ${\rm Gl}_p(K)$ making $D_yH(0,0)$ unipotent upper triangular. In the sequel we shall always assume that $\inin(H_i)=y_i$ with respect to the chosen monomial order on $\N^{n+p}$. \med

The {\it baby series vector} of a mother code $H\in K[x,y]^p$ is the formal power
series vector $h=(h_1\to h_p)\in \kxx^p$ vanishing at $0$ which is the unique solution
of $H(x,h(x))=0$. The existence and uniqueness of $h$ are ensured by the implicit function theorem for formal power series. The components $h_i$ of the baby series vector $h$ are algebraic series. This can be seen by the algebraic implicit function theorem [KPR, p.÷ 91], or Artin's Approximation Theorem [Ar2], or by the following argument: Consider the system $H(x,y)=0$ as equations for the last variable $y_p$. After a renumeration of the components of $H$, the last derivative $\6_{y_p}H_p(0,0)$ does not vanish. There then exists a unique solution $h_p(x,y_1\to y_{p-1})$ of $H_p(x,y_1\to y_{p-1},y_p)=0$ vanishing at $0$, and $h_p$ is algebraic over $K[x,y_1\to y_{p-1}]$. By induction on $n$ and the transitivity of algebraicity we conclude that $h=(h_1\to h_p)$ is algebraic. \med

%-------------------------------------------------
%     BEGIN  TEXT  MARIEMI 
%-------------------------------------------------
\ignore

[Include here the direct argument of ME: Es siempre la historia de que cuando eliminas algunas 
variables en un sistema polinomial $ H_j(x_1,...,x_n,y_1,...,y_p)=0$, $j=1,...,p 
$ con la condicion jacobiana no nula. Por ejemplo si eliminamos las 
$y_2,...,y_p$ estoy 
considerando la  proyeccion de un conjunto algebraico \med

\hs 3cm $V= \{H_j=0  \}\subset K^{n+p}$\med

sobre $K^{p+1}$ que tiene una componente 
$W$ irreducible en $0$ de dim $n$ (por la condicion jacobiana). Cerca del 
$(0,...,0) \in K ^n$ la proyeccion de $W$ no puede tener mayor dimension de 
$n$. Asi como $p\geq n$ es $p+1\geq n$ y existe una ecuacion algebraica no 
trivial  $P(x_1,...,x_n,y_1)=0$ en la imagen.] \med

\recognize
%-------------------------------------------------
%    END  TEXT  MARIEMI 
%-------------------------------------------------

An algebraic power series vector $h=(h_1\to
h_p)\in \kxx^p$ is {\it a baby series vector} if it admits a mother code $H\in
K[x,y]^p$ defining it.\med

A {\it father code} is a vector $G=(G_1\to G_r)$ of polynomial vectors $G_i\in
K[x,y]^s$ (there are no further conditions on the $G_i$). We consider $G$ as a row
vector with entries the column vectors $G_i$, say as a matrix in 
$K[x,y]^{s\times r}$. \med

A {\it family code} is a pair $(H,G)$ where $H\in K[x,y]^p$ is a
mother code and $G\in K[x,y]^{s\times r}$ a father code, both carrying on the
same sets of variables. We say that algebraic power series vectors $g_1\to g_r\in\kxx^s$ {\it have family code}
$(H,G)\in K[x,y]^p\times  K[x,y]^{s\times r}$ if \med

\hs 4cm $g_k=G_k(x,h(x))$ \med

for $1\leq k\leq r$, where $h\in \kxx^p$ is the baby series vector of the mother code
$H$. The vectors $g_k$ hence belong to $K[x,h]^s\subset K[[x]]^s$. We call $h$ the baby
series vector underlying $g_1\to g_r$, or, the other way round, $g_1\to g_r$ the
algebraic power series vectors produced from $h$ by the father code $G$.\med

%-----------------------------------------------

{\bf Example 6.1.} Take $g_1=z^3+z^2h$, $g_2=xz^2+xzh$ with baby series
$h=1-\sqrt{1-x^2}$, mother code $H=2y-y^2-x^2$ and father code
$G_1=z^3+z^2y$,  $G_2=xz^2+xzy$. Notice that the second series solution
$1+\sqrt{1-x^2}$ of $H=0$ has non-zero constant term and is therefore not
considered as a baby series of $H$.\med

%-------------------------------------------------

{\bf Example 6.2.} Let $H$ be the vector  $(H_1,H_2)$ with
$H_1=y_1^2-2y_1-y_2-x_2$ and $H_2=x_2y_2^2-2y_2-y_1-x_1$. The
vector $H$  is the mother code of the baby series vector $(h_1,
h_2)$ where $h_1$ and $h_2$ are related by $h_1=1-\sqrt{1+x_2+h_2}$ and
$h_2=(1-\sqrt{1+x_2(x_1+h_1)})/x_2$. The mother code $H$ defines the same baby series vector as  the mother code $H'=(H'_1,H'_2)$ given by \med

\hs 1cm $H'_1 =
-x_1+x_2^3+2x_2-4x_2y_1^3+(3+4x_2^2)y_1+(-2x_2^2-2+4x_2)y_1^2+x_2y_1^4$,\med

\hs 1cm  $H'_2=  -y_1^2+2y_1+y_2+x_2$. \med

Now, $D_yH'(0,0)$ is unipotent upper triangular and $H_1'$ does not depend on $y_2$. Hence, the expansion
of the series $h_1$ can be computed up to a any order from the equation $H_1'=0$. From
$H_2'=0$ we get $h_2=-x_2-2h_1+h_1^2$.\med

\big\goodbreak

%-----------------------------------------------
%           CONSTRUCTION OF CODES
%-----------------------------------------------

{\Bf 7. Construction of codes}\med

Codes of algebraic power series as above were introduced by Alonso, Mora and
Raimondo. Their construction is based on an effective version of the Artin-Mazur theorem [AM, p.÷ 88, AMR, appendix, BCR, Thm.÷ 8.4.4, p.÷ 173].  \big

{\bf Theorem 7.1.} {\it For any algebraic series $g\in \kxx$ there is a finite algorithm to construct from an algebraic relation $P(x,t)=0$ satisfied by $g$ and the Taylor expansion of $g$ up to sufficiently high degree a family code $(H,G)\in
K[x,y]^p\times K[x,y]$ of $g$, for some $p$.}\big

%-----------------------------------------------

{\it Proof.} Let $P(x,g(x))=0$ be a minimal hence irreducible algebraic relation for $g$. Denote by $X\sub \A_K^{n+1}$ the zero-set of $P$ in affine $(n+1)$-space $\A_K^{n+1}$ over $K$. We assume that $g(0)=0$ so that $(0,0)\in X$. Let $Y$ be the normalization of $X$. Choose an embedding $Y\subset \A^{n+p}$ so that the normalization map $\pi:Y\map X$ is induced by the projection $\A^{n+p}\map \A^{n+1}$, $(x,y)\map (x,y_1)$ on the first $n+1$ components.\med

The Taylor expansion of $g$ specifies a unique point $b\in Y$ which maps to $0\in X$ and through which, by the universal property of normalization, passes a lifting $(x,\widetilde g(x))$ of $(x,g(x))$. From Zariski's Main Theorem [Za, Mu, p.÷ 209] we know that $Y$ is analytically irrreducible at $b$. But as $Y$ contains the graph of $\widetilde g$ and has dimension $n$, it is smooth at $b$. By the Jacobian criterion it is therefore possible to choose polynomial equations $H_1\to H_p$ defining $Y$ in a Zariski neighborhood of $b$ in $\A^{n+p}$ and satisfying at $b$ the assumption of the implicit function theorem, i.e., of a mother code. Let  $(h_1\to h_p)$ be the associated baby series vector. By the special choice of $\pi$ we get $g=h_1$, say $g=G(h_1\to h_p)$ with father code $G=y_1$. This proves the theorem.\med

The construction of the normalization is effective [dJP] and implemented for instance in the computer-algebra program Singular [GPS].\med

When handling several algebraic power series it is more economic to work with one mother code and several father codes instead of choosing separate mother codes for each series. This goes as follows.\med

Let be given mother codes $H^{j}\in K[x,y^{j}]^{p_j}$ for $j=1\to r$ in distinct
sets of variables $y^{j}=(y_1^{j}\to y_{p_j}^{j})$ defining baby series vector
$h^{j}=(h_1^{j}\to h_{p_j}^{j})\in \kxx^{p_j}$. The {\it direct sum} $H$ of
the $H^{j}$'s is given as the row vector $H=(H^1\to H^r)\in \Pi_{j=1}^r
K[x,y]^{p_j}\isom K[x,y]^p$, where $y$ denotes the collection of all $y^{j}$ and
$p=\sum p_j$. This $H$ is again a mother code, because the Jacobian matrix
$D_yH(0,0)$ of $H$ with respect to $y$ at $0$ has block diagonal form with
invertible blocks equal to $D_{y^{j}}H^{j}(0,0)$ on the diagonal. The vector
$h=(h^{1}\to h^{r})$ obtained by listing all baby series vectors $h^j$ of the
mother codes $H^{j}$ in a row is the baby series vector of $H$. This passage to
direct sums of mother codes allows us to treat several baby series vectors $h^{j}$
simultaneously as one baby series vector $h$ (with many components).
Accordingly, finitely many algebraic series can always be considered as
produced by certain father codes from the {\it same} baby series vector
$h=(h_1\to h_p)$ of {\it one} mother code $H\in K[x,y]^p$. This allows us to work
throughout with vectors in $K[x,h_1\to h_p]^s$. \med

Note that mother codes as defined above may require large sets of variables and are thus computationally very expensive.\med

%[Compare with remarks of ME and with Iversen, LNM 310, Cor. 2.2.]\med

\big\goodbreak

%-----------------------------------------------
%          CODES FOR MODULES
%-----------------------------------------------

{\Bf 8. Codes for modules of algebraic series}\med

Let be given algebraic power series vectors $g_1\to g_r\in \kxx^s$ vanishing at $0$ with mother code $H\in K[x,y]^p$, baby series vector $h\in \kxx^p$ and father code $G\in
K[x,y]^{s\times r}$ so that $g_k=G_k(x,h(x))$. The submodule $\<g_k\>$ of $K[[x]]^s$ generated by the series $g_k$ admits the following polynomial description. \med

%-----------------------------------------------

{\bf Lemma 8.1.} {\it Let $\<(y_i-h_i)\cdot e_\ell,g_k\>$ and $\<H_i\cdot e_\ell,
G_k\>$ be the submodules of $K[[x,y]]^s$ generated by the vectors $(y_i-h_i)\cdot e_\ell$
and $g_k$, respectively $H_i\cdot e_\ell$ and $G_k$, for $1\leq i\leq p$, $1\leq
\ell\leq s$, $1\leq k\leq r$. Then} \med

\hs 3cm $\<(y_i-h_i)\cdot e_\ell,g_k\>  \,=\,  \<H_i\cdot e_\ell, G_k\>$.\med

%-----------------------------------------------

{\it Proof.} We fix a monomial order $<_\eta$ on $\N^n\times
\{1\to s\}$ and choose an extension  $<_\eps$ of $<_\eta$ to $\N^{n+p}\times
\{1\to s\}$ which is degree compatible with respect to $\N^p$ and satisfies
$y_i \cd e_\ell<_\eps x_j\cd e_\ell$ for all $1\leq i\leq p$, $1\leq j\leq n$ and $1\leq
\ell\leq s$. After a suitable multiplication of $H$ with a constant matrix in ${\rm GL}_p(K)$ we may assume that $\inin(H_i\cd e_\ell) = y_i\cd e_\ell$.\med

The ideal $\<H_i\>$ of $K[[x,y]]$ generated by $H_1\to H_p$ is  contained in the
ideal $\<y_i-h_i\>$ because of $H(x,h(x))=0$. Take a monomial order $<_\delta$ on $\N^{n+p}$ so that $y_i<_\delta x_j$ for all $1\leq i\leq p$ and $1\leq j\leq n$.
The initial ideals of $\<H_i\>$ and $\<y_i-h_i\>$ coincide because, by the choice of $<_\delta$, they are both generated by $y_1\to y_p$. By the Division Theorem for formal power series, the two ideals coincide. As $g_k$ is obtained from $G_k$ by replacing $y_i$ by $h_i$, the submodules of
$K[[x,y]]^s$ generated by $g_1\to g_r$, respectively $G_1\to G_r$ are congruent
modulo $\<y_i-h_i\>\,=\,\<H_i\>$. This proves the lemma. \med

%-------------------------------------------------

%\[Comment (18) from may 30 2009: We even have $\<y_i-h_i\>=\<H_i\>$ as ideals in $K[[x,y]]$. The proof should not need the division theorem (with respect to a suitably chosen order on $\N^{n+p}$). Instead, the equality should follow directly from the Inverse Function Theorem. From the equality $\<y_i-h_i\>=\<H_i\>$ then follows immediately the equality $\<y_i-h_i, g_k\>=\<H_i, G_k\>$ of the lemma, because the $G_k$ reduce to $g_k$ modulo substitution of $y_i$ by $h_i$.]\med

%-------------------------------------------------

We call $\widetilde I=\,\<H_i\cdot e_\ell, G_k\>\, \subset K[[x,y]]^s$, or, more accurately, its polynomial generators $H_i\cdot e_\ell$ and $G_k$, the {\it family code} of the submodule  $I=\,\<g_k\>$ of $\kxx^s$. Observe that $\widetilde
I\cap \kxx^s=I$. \med

%-----------------------------------------------

{\bf Lemma 8.2.} {\it Let be given a monomial order $<_\eta$ on $\N^n\times
\{1\to s\}$ and an extension  $<_\eps$ of $<_\eta$ to $\N^{n+p}\times
\{1\to s\}$ which is degree compatible with respect to $\N^p$ and satisfies
$ y_i \cd e_\ell<_\eps x_j\cd e_\ell$ for all $1\leq i\leq p$, $1\leq j\leq n$ and
$1\leq \ell\leq s$. Let $\widetilde I=\,\<H_i\cdot e_\ell, G_k\>$  and
$I=\,\<g_k\>$ be the respective submodules of $K[[x,y]]^s$ and $\kxx^s$.
Then} \med

\hs 3cm $\inin(\widetilde I) \cap \kxx^s=\inin(I)$.\med

%-----------------------------------------------

{\it Proof.} We may choose a minimal reduced standard basis of $\widetilde I$. Let $\widetilde g_k$ be an element of this basis which does not have an initial monomial vector of the form $y_i\cd e_\ell$. From reducedness it follows that $\widetilde g_k$ is independent of $y_1\to y_p$, say $\widetilde g_k\in \widetilde I \cap \kxx^s=I$. In particular, the vectors $\widetilde g_k$ form a standard basis of $I$ and hence $\inin(\widetilde I) \cap \kxx^s=\inin(I)$. \med

\big\goodbreak

%-----------------------------------------------
%     CONSTRUCTION OF STANDARD BASIS
%-----------------------------------------------

{\Bf 9. Construction of standard basis}\med

The first construction we need is a direct consequence of Mora's tangent cone algorithm
[Mo], respectively Lazard's homogenization method [Lz],  cf.÷ also  with
[AMR, Thm.÷ 1.3, CLO, p.÷ 202, Gr1, Gr2, GP, Thm.÷ 6.4.3]. It provides an algorithm to construct the family code of a (not necessarily reduced) standard basis of a module of algebraic power series vectors. \med

%-----------------------------------------------

{\bf Theorem 9.1.} {\it Let $I$ be a submodule of $\kxx^s$ generated by
algebraic power series vectors $g_1\to g_r\in \kxx^s$ which are given by their family code.
Let be chosen a monomial order $<_\eta$ on $\N^n\times \{1\to s\}$.
There is a finite algorithm to compute the family codes of the elements of a standard
basis of $I$ with respect to $<_\eta$ from the family codes of $g_1\to g_r$.
In particular, it is possible to compute the initial module $\inin(I)$ of
$I$.}\med

%-----------------------------------------------

{\it Proof.} Let $g_1\to g_r$ have mother code $H\in K[x,y]^p$, baby series
vector $h\in \kxx^p$ and father code $G\in K[x,y]^{s\times r}$. Extend $<_\eta$
to a monomial order $<_\eps$ on $\N^{n+p}\times \{1\to s\}$ which is degree
compatible with respect to $\N^p$ and satisfies $y_i\cd e_\ell<_\eps
x_j\cd e_\ell$ for all $i$, $j$ and $\ell$. We assume w.l.o.g. that the
initial monomial vectors of $H_i\cd e_\ell$ with respect to $<_\eps$ are $y_i\cd e_\ell$. \med

As  $\widetilde I=\, \<H_i\cdot e_\ell, G_k\>$ is generated by polynomial vectors,
Mora's  tangent cone algorithm or Lazard's homogenization method apply to
construct a polynomial standard basis for it. This basis is in general not reduced. We may choose a minimal basis consisting of the vectors $H_i\cd e_\ell$ with $\inin(H_i\cd e_\ell)= y_i\cd e_\ell$ and of other polynomial vectors $\widetilde G_1\to \widetilde G_{r'}\in K[x,y]^s$ with initial monomial vectors in $\kxx^s$. The latter form the father code of algebraic power series vectors $\widetilde g_1\to \widetilde g_{r'}\in\kxx^s$, say $\widetilde g_k=\widetilde G_k(x,h)$. Note that $\widetilde G_k$ is congruent to $g_k$ modulo the submodule $\<H_i\cd e_\ell\>$ of $K[[x]]^s$. By Lemma 8.2, the $\widetilde g_k$ form a standard basis of $I$. This proves the theorem. \med

\big\goodbreak

%-----------------------------------------------
%     CONSTRUCTION OF REDUCED STANDARD BASIS
%-----------------------------------------------

{\Bf 10. Construction of reduced standard basis}\med

The central part in establishing the division algorithm for modules of
algebraic power series vectors is the construction of a {\it reduced} standard basis. The mere existence follows from Hironaka's theorem. The effective part in the special case of principal ideals, i.e., the construction of the code of the Weierstrass form of an $x_n$-regular algebraic power series, has been established by Alonso, Mora and Raimondo [AMR, Thm.÷ 5.5]. The general statement is as follows:\med

{\bf Theorem 10.1.} {\it Let $I$ be a submodule of $\kxx^s$ generated by
algebraic power series vectors. Assume that $I$ satisfies Hironaka's box 
condition with respect to a monomial order $<_\eta$ on $\N^n\times \{1\to s\}$.
Then the family codes of a reduced standard basis of $I$ can be computed by a
finite algorithm from the family codes of any algebraic power series vectors
$g_1\to g_r\in \kxx^s$ generating $I$.}\med

%-----------------------------------------------

The proof of this result is given in sections 13 to 15. In the formal power series case, a reduced standard basis can be constructed up to any given degree by dividing monomial generators of the initial module by the module itself. For algebraic series, this construction would require to dispose already of an effective division algorithm. To avoid this logical cycle, reduced standard bases have to be constructed in a different way. \med

The clue relies in the concept of a {\it virtual reduced standard basis}. 
Such a basis consists of polynomial vectors whose coefficients are unknown and written as new variables. Upon replacing the variables by suitable series
in $x$, the virtual reduced standard basis will transform into an actual reduced standard
basis of the module. The resulting coefficient series of the actual reduced standard
basis -- more precisely, their codes -- are computed by dividing the polynomial generators of the module $\<(y_i-h_i)\cdot e_\ell,g_k\>\,=\,\<H_i\cdot e_\ell, G_k\>$  by the virtual
basis using the polynomial division algorithm. The definition requires that both the generators and the virtual basis are polynomial vectors, and that the initial monomial vectors of the virtual reduced standard basis can be interpreted as the leading (i.e., maximal) monomial vectors w.r.t.÷ another, suitably chosen monomial order. The choice of this order is rather subtle, see section 13. The remainders of the division then allow us to extract the codes of the required coefficients series.\med

\big\goodbreak

%-----------------------------------------------
%     EFFECTIVE  DIVISION  FOR ALGEBRAIC POWER SERIES
%-----------------------------------------------

{\Bf 11. Effective division for algebraic power series}\med

Our main result asserts that the division by modules of algebraic power series vectors
with box condition can be made effective, i.e., can be performed by applying finitely many operations to the codes. The case of principal ideals $I$, say the effective Weierstrass Division Theorem for algebraic power series, is due to Alonso, Mora and Raimondo in [AMR, Thm.÷ 5.6].\med

%-----------------------------------------------

{\bf Theorem 11.1.} {\it Let $I$ be a submodule of $\kxx^s$ generated by algebraic
power series vectors. Assume that $I$ satisfies Hironaka's box condition with
respect to a monomial order $<_\eta$ on $\N^n\times \{1\to s\}$. Let be given
the family codes of algebraic power series vectors $g_1\to g_r\in\kxx^s$
generating $I$. There exists a finite algorithm which computes, for any algebraic power series vector $f\in\kxx^s$, from the family code of $f$ the family codes of
algebraic power series $a_1\to a_r$ in $\kxx$ and of an algebraic power series
vector $c\in\co(I)\subset\kxx^s$ so that \med

\hs 4cm $f= \sum_{k=1}^r a_kg_k +c$\med

is a formal power series division of $f$ by $g_1\to g_r$.}\med

%-----------------------------------------------

The algorithm produces quotients $a_k$ which are algebraic series but which in
general need not satisfy the support conditions of the formal power series
division as in [Ga]. The remainder $c$, of course, is unique and only
depends on the chosen monomial order $<_\eta$. \med

We shall prove Theorem 11.1 by first constructing from $g_1\to g_r$ via Theorems 9.1 and 10.1 the family codes of a reduced standard basis of $I$. The division
algorithm for a reduced standard basis will then be established by
induction on the number of variables.\med

\big\goodbreak
%-------------------------------------------------
%          PROOFS
%-------------------------------------------------

%\cl {\Bf PROOFS}\big

%-----------------------------------------------
%           OUTLINE PROOF OF THM 2  AND 3  
%-----------------------------------------------

{\Bf 12. Logical structure of the proofs of Theorems 10.1 and 11.1}\med

Both theorems will be established independently of Hironaka's existential
division theorem. We start with establishing Theorem 10.1, the construction of the codes of a reduced standard basis, in the special case of $x_n$-regular modules. This is the hardest part of the whole story. It relies on introducing the virtual reduced standard basis of the module, which allows us to perform {\it polynomial} divisions for constructing the required codes. This section is inspired by Mora's tangent cone algorithm and the techniques of Alonso, Mora and Raimondo in [AMR]. Extracting from the virtual reduced standard basis the actual reduced standard basis uses  in an essential way the assumption of $x_n$-regularity.\med

From Theorem 10.1 for $x_n$-regular modules we deduce the division algorithm of Theorem 11.1 for $x_n$-regular modules. The algorithm uses again a polynomial
division, this time by the codes of the reduced standard basis. For its
termination it is necessary that the basis is already reduced. \med

The general cases of Theorems 10.1 and 11.1 are then deduced simultaneously from the special cases by induction on the number of variables and using Hironaka's box condition  together with the notion of Janet basis. One selects from the given (not yet reduced) standard basis of the module those elements which are $x_n$-regular. Such elements exist because of the box condition. Considering the module generated by these elements, one may construct the codes of its reduced standard basis via Theorem 10.1 in the special case. Then Theorem 11.1 allows us to reduce effectively the remaining elements with respect to the first set of elements. By induction on the number of variables, the tails of the first elements can now be divided conversely by the remaining elements, yielding eventually the codes of the whole reduced standard basis of the module. Once this is achieved, it is relatively simple to establish also the division of Theorem 11.1 in the general case. \med 

\big\goodbreak

%-----------------------------------------------
%   PROOF OF THM 10.1  FOR $x_n$-REGULAR MODULES
%-----------------------------------------------

{\Bf 13. Proof of Theorem 10.1 for $x_n$-regular modules}\med

In the situation of Theorem 10.1, we first treat  the case where $I$ is
$x_n$-regular with respect to $<_\eta$. As seen in Lemmata 8.1 and 8.2, it suffices to construct the family code of a reduced standard basis of the submodule $\widetilde I=\,\<H_i\cdot e_\ell, G_k\>\,=\,\<(y_i-h_i)\cdot e_\ell,g_k\>$ of $K[[x,y]]^s$ with respect to the chosen extension $<_\eps$ of $<_\eta$. Note here that $\widetilde I$ is not $x_n$-regular, since also the $y_i\cd e_\ell$ appear in the initial module. This is, however, not a serious drawback. The proof is somewhat involved and goes in several steps. Let us first specify the setting.\med

(a) We suppose that the generators $g_k$ of $I$ vanish at $0$ for all $1\leq k\leq r$. Hence this also holds for all $H_i\cd e_\ell$ and $G_k$. We may assume by Theorem 9.1 and its proof that the polynomial vectors $H_i\cdot e_\ell$
and $G_k$ form a minimal standard basis of $\widetilde I$. As $I$ is
$x_n$-regular  and $\inin(H_i\cdot e_\ell)=y_i\cdot e_\ell$, the initial module
$\inin(\widetilde I)$ is generated by $y_i\cdot e_\ell$ and monomial vectors
$x_n^{d_k}\cd e_{m_k}$ for some $d_k> 0$, $1\leq m_k\leq s$ and $1\leq k\leq r$. Hence $r\leq s$. After a suitable permutation of the components of $\kxx^s$ and a renumeration of $G_1\to G_r$ we may assume that $m_k=k$, say $\inin(G_k)=x_n^{d_k}\cd e_k$ for all $k$. The permutation of the components is only for notational convenience. It will not affect the induction we shall apply later on when proving Theorems 10.1 and 11.1 in the general case.\med

The canonical direct monomial complement $\co(\widetilde I)$ of $\inin(\widetilde I)$ in
$K[[x,y]]^s$ is of form\med

\hs 2cm $\co(\widetilde I)=\oplus_{m=1}^r\,\oplus_{j=0}^{d_m-1} K[[x']]\cdot x_n^j\cd e_m \,\, \oplus\,\,\oplus_{m=r+1}^s\, K[[x]]\cd e_m$.\med

Write a reduced standard basis of $\widetilde I$ as \med 

\hs .7cm $b_{i\ell}= y_i\cdot e_\ell - b_{i\ell}^\circ -
\sum_{m=1}^r\sum_{j=0}^{d_m-1}  u_{i\ell mj}(x')\cdot x_n^j\cdot e_m -
\sum_{m=r+1}^s  v_{i\ell m}(x)\cdot e_m$, \med 

\hs .7cm $b_k= x_n^{d_k}\cdot e_k - b_k^\circ -\sum_{m=1}^r\sum_{j=0}^{d_m-1}
 u_{kmj}(x')\cdot x_n^j\cdot e_m - \sum_{m=r+1}^s  v_{km}(x)\cdot e_m$,
\med

with polynomial vectors $b_{i\ell}^\circ$ and $b_k^\circ$ in\med

\hs 3cm 
$\oplus_{m=1}^r\,\oplus_{j=0}^{d_m-1} K x_n^j\cd e_m \,\, \oplus\,\,
\oplus_{m=r+1}^s\, K\cd e_m$ \med

and power series $u_{i\ell mj}(x')$, $v_{i\ell m}(x)$, $u_{kmj}(x')$ and $v_{km}(x)$ vanishing at $0$. It is necessary here to split off $b_{i\ell}^\circ$ and $b_k^\circ$ because the mother codes of algebraic power series are only defined for series vanishing at $0$. Note that $u_{i\ell mj}(x')$ and $u_{kmj}(x')$ do not depend on $x_n$, and that $b_{i\ell}^\circ$ and $b_k^\circ$ vanish at $0$ because the $H_i\cd e_\ell$ and $G_k$ do. In particular, these vectors have zero entries in the last $s-r$ components. Since $\inin_\varepsilon(b_{i\ell})=y_i\cd e_\ell$ and $\inin_\varepsilon(b_k)=x_n^{d_k}\cd e_k$ the $\ell$-th component of $b_{i\ell}^\circ$ and the $k$-th component of $b_k^\circ$ are both zero.\med

The series $b_{i\ell}$ have different shapes according to whether $1\leq \ell\leq r$ or $r+1\leq\ell\leq s$. Namely, for $r+1\leq\ell\leq s$, it follows from the $x_n$-regularity of $I$ that the vectors $(y_i-h_i)\cd e_\ell$ are already reduced. Hence we have  $b_{i\ell}=(y_i-h_i)\cd e_\ell$ for $r+1\leq\ell\leq s$. This will be used later on. We are grateful to D.÷ Wagner for specifying an inaccuracy which appeared at this place in an earlier draft of the paper.\med

In a first step, we determine the vectors $b_{i\ell}^\circ$ and $b_k^\circ$. Afterwards, the family codes of the coefficient series $u_{i\ell mj}(x')$, $v_{i\ell m}(x)$, $u_{kmj}(x')$ and $v_{km}(x)$ will be constructed. This will show in
particular that they are algebraic series. \med

%-----------------------------------------------

(b) In order to compute $b_{i\ell}^\circ$ and $b_k^\circ$, one can
construct the reduced standard basis of $\widetilde I$ up to a sufficiently high degree by
applying its formal power series construction modulo a sufficiently high power
of the maximal ideal of $K[[x,y]]$.
\med

%-----------------------------------------------

(c) The series $u_{i\ell mj}(x')$, $v_{i\ell m}(x)$, $u_{kmj}(x')$
and $v_{km}(x)$ will be determined by a trick which has already appeared in the literature, see e.g.÷ [AMR]: Define the {\it virtual reduced standard basis} of $\widetilde I$ as the
polynomial vectors\med

\hs 1cm $B_{i\ell}= y_i\cdot e_\ell - b_{i\ell}^\circ-
\sum_{m=1}^r\sum_{j=0}^{d_m-1} u_{i\ell mj}\cdot x_n^j\cdot e_m - \sum_{m=r+1}^s
v_{i\ell m}\cdot e_m$, \med 

\hs 1cm $B_k= x_n^{d_k}\cdot e_k - b_k^\circ - \sum_{m=1}^r\sum_{j=0}^{d_m-1}
u_{kmj}\cdot x_n^j\cdot e_m - \sum_{m=r+1}^s v_{km}\cdot e_m$, \med

where $u_{i\ell mj}$, $v_{i\ell m}$, $u_{kmj}$ and $v_{km}$ are now new
variables (to be abbreviated by $u$ and $v$). From these we shall construct
certain polynomials $U_{i\ell mj}$, $V_{i\ell m}$, $U_{kmj}$ and $V_{km}$ in
$K[x,y,u,v]$. All these together will constitute the mother code $(U,V)$  of
the  baby series vector $(u(x'),v(x))$ of components $u_{i\ell
mj}(x')$,  $u_{kmj}(x')$, respectively  $v_{i\ell m}(x)$, $v_{km}(x)$.  And, 
consequently, $B_{i\ell}$ and $B_k$ will be the father codes of the series
vectors $b_{i\ell}$ and $b_k$ we were looking for, with $b_{i\ell},\, b_k\in
K[x,y,u(x'),v(x)]^s$.\med

We have noticed above that, for $r+1\leq\ell\leq s$, the vectors $b_{i\ell}$ equal $(y_i-h_i)\cd e_\ell$. As the polynomial vectors $B_{i\ell}$ are the father codes of $b_{i\ell}$ they will therefore have, for $r+1\leq\ell\leq s$, only one non-zero entry, namely in the $\ell$'s component. Hence we may set all variables $u_{i\ell mj}$, $v_{i\ell m}$ for $r+1\leq m\leq s$ and $m\neq\ell$ equal to $0$. This will be used below when proving the independence of $U_{i\ell mj}$ and $U_{kmj}$ on $v$.\med

%-----------------------------------------------

(d) The construction of the codes $U$ and $V$ uses the polynomial division
algorithm with respect to a suitably chosen monomial order. Compare monomial vectors $u^\gamma v^\delta x^\alpha y^\beta\cd e_m$ by considering the integer vector \med

\hs 3cm $\ds (\beta, \alpha_n- d_m , \alpha', -m, \gamma, \delta)$ \med

lexicographically. Here, the tuples $\gamma$ and $\delta$ are taken as ordered vectors, e.g.÷ by choosing some ordering of their indices. It is easily checked that this defines a monomial order $<_\omega$ on $\N^{q+n+p}$, where $q$ is the number of $u$ and $v$ variables, and that the leading (= maximal) monomial vectors of $B_{i\ell}$ and $B_k$ w.r.t.÷ $<_\omega$ are $y_i\cd e_\ell$ and $x_n^{d_k}\cdot e_k$.\med

We now divide $H_i\cdot e_\ell$ and $G_k$ polynomially by $B_{\i\l}$ and $B_\k$ with respect to this monomial order, say with leading monomial vectors
$y_\i\cdot e_\l$ and $x_n^{d_\k}\cdot e_\k$ and the scopes $n_{\i\l}=q+n+\i$ and $n_\k=q+n$ (with $1\leq \i\leq p$, $1\leq \l\leq s$ and $1\leq \k\leq r$). The division
yields in finitely many steps remainders $R_{i\ell}$ and $R_k$ in the canonical
direct monomial complement \med

\hs .5cm $K[u,v]\otimes\co(\widetilde I) =\oplus_{m=1}^r(\oplus_{j=0}^{d_m-1}
K[u,v][[x']]\cdot x_n^j)\cdot e_m \,\, \oplus \,\, \oplus_{m=r+1}^s K[u,v][[x]]\cdot e_m$ \med

of $K[u,v]\otimes \inin(\widetilde I)$ in $K[u,v][[x,y]]^s$. Expanding these
remainders as polynomial vectors in $x_n$ yields\med

\hs 2cm $R_{i\ell}=\sum_{m=1}^r \sum_{j=0}^{d_m-1} U_{i\ell m j}\cdot x_n^j\cdot
e_m + \sum_{m=r+1}^s V_{i\ell m}\cdot e_m$,\med

\hs 2cm $R_k=\sum_{m=1}^r \sum_{j=0}^{d_m-1} U_{k m j}\cdot x_n^j\cdot e_m +
\sum_{m=r+1}^s V_{k m}\cdot e_m$, \med

with polynomials $U_{i\ell mj}$,  $U_{kmj}$ in $K[u,x']$ and  $V_{i\ell
m}$, $V_{km}$ in $K[u,v,x]$. Note here that $U_{i\ell mj}$ and $U_{kmj}$ do
not depend on $v$ because $v_{i\ell m}$ and $v_{km}$ only appear in the last
$s-r$ components of $B_{i\ell}$ and $B_k$ and because $u_{i\ell mj}$ and $v_{i\ell m}$ can a priori be set equal to $0$ for $m\neq\ell$.\med

%-----------------------------------------------

(e) We show that $U$ and $V$ have no constant terms. Replacing in $R_{i\ell}$
and $R_k$ the variables $u$ and $v$ by the series $u(x')$ and $v(x)$ produces
power series vectors $r_{i\ell}$ and $r_k$ which belong to $\co(\widetilde I)$
because $u(x')$ does not depend on $x_n$ and $U$ does not depend on $v$. But
by construction, $r_{i\ell}$ and $r_k$ also belong to $\widetilde I$. From the
formal power series division follows that both $r_{i\ell}$ and $r_k$ are
identically zero. This in turn implies by the direct sum decomposition of $\co(\widetilde I)$ that replacing in $U$ and $V$  the
variables $u$ and $v$ by $u(x')$ and $v(x)$ gives zero. As $u(x')$ and $v(x)$
have no constant term, also $U$ and  $V$ have no constant term.\med
 
%-----------------------------------------------
%       START     NEW  
%-----------------------------------------------

(f) We show that $U$ and $V$ form a mother code of certain baby series. By the description of mother codes it suffices to find a monomial order $<_\xi$ on $\N^{q+n+p}\times\{1\to s\}$ such that the respective initial monomials of $U_{i\ell mj}(0,0,u,v)$, $V_{i\ell m}(0,0,u,v)$, $U_{kmj}(0,0,u,v)$ and 
$V_{km}(0,0,u,v)$ are $u_{i\ell mj}$, $v_{i\ell m}$,
$u_{kmj}$ and $v_{km}$. By taking an order which is compatible with the degree in the $u$ and $v$ variables it suffices to prove the above for the linear parts of $U_{i\ell mj}(0,0,u,v)$, $V_{i\ell m}(0,0,u,v)$, $U_{kmj}(0,0,u,v)$ and 
$V_{km}(0,0,u,v)$. \med

These linear parts are given by the first substitution step of the polynomial division as the coefficients of  $x_n^j\cd e_m$ (with $1\leq j\leq d_m-1$, $1\leq m\leq r$) respectively $e_m$ (with $r+1\leq m \leq s$), when dividing $H_i\cd e_\ell$ and $G_k$ by the vectors $B_{\iota\lambda}$ and $B_\kappa$ ($1\leq \iota \leq p$, $1\leq \lambda \leq s$, $1\leq \kappa \leq r$) with leading monomial vectors $y_\iota\cd e_\lambda$ and $x_n^{d_\kappa}\cd e_\kappa$ and scopes $q+n+\iota$, respectively $q+n$. Here, the $y$ variables are ordered naturally $y_1\to y_p$, so that the scope $q+n+\iota$ of $y_\iota\cd e_\lambda$  allows multiplication of $B_{\i\l}$ with polynomials in $x_1\to x_n$, $y_1\to y_\i$ and all $u$ and $v$ variables. \med

Note that the polynomial vectors $b_{\i\l}^\circ$ and $b_\k^\circ$ of $K[x]^s$ appearing in $B_{\i\l}$ and $B_\k$ vanish at zero and hence do not contribute to the linear terms of $U_{i\ell mj}(0,0,u,v)$, $V_{i\ell m}(0,0,u,v)$, $U_{kmj}(0,0,u,v)$ and 
$V_{km}(0,0,u,v)$.\med

The construction of the monomial order $<_\xi$ on $\N^{q+n+p}\times\{1\to s\}$ involves a monomial order $<_\zeta$ on $\N^q$ (recall that $q$ is the number of $u$ and $v$ variables) whose choice is motivated by the following computations (where we shall assume throughout w.l.o.g.÷ that $\inin(H_i\cd e_\ell)=y_i\cd e_\ell$ and $\inin(G_k)=x_n^{d_k}\cd e_k$).\med

%-------------------------------------------------

{\it Linear terms of $\, U_{i\ell mj}(0,0,u,v)$}: These occur after the first substitution step of the polynomial division as the coefficients of  $x_n^j\cd e_m$ (with $1\leq i\leq p$, $1\leq \ell\leq s$, $1\leq m\leq r$, $1\leq j\leq d_m-1$) when dividing $H_i\cd e_\ell$  by the vectors $B_{\i\l}$ and $B_\k$ with leading monomial vectors $y_\i\cd e_\l$ and $x_n^{d_\k}\cd e_\k$ and scopes $q+n+\i $, respectively $q+n $  (where $\i $, $\l$ and $\k$ vary in the ranges  $1\leq \i \leq p$, $1\leq \l \leq s$, $1\leq \k \leq r$). Notice that the polynomials $U_{i\ell mj}$ do not depend on $y$ and $v$.\med

Let $x^\rho y^\sigma\cd e_\ell$ be a monomial vector of the expansion of $H_i\cd e_\ell$, with  $\rho\in \N^n$, $\sigma\in \N^p$. If it is a multiple of the leading monomial vectors $y_\i\cd e_\l$, respectively $x_n^{d_\k}\cd e_\k$, of $B_{\i\l}$, respectively $B_\kappa$, subject to the correct scope conditions, it will be replaced in the polynomial division by the according multiple of the tails $\ol B_{\i\l}$, respectively $\ol  B_\k$. After the substitution we have to look at the coefficient of $x_n^j\cd e_m$ and set $x=0$ and $y=0$. We distinguish three cases.\med

(i) The substitution of the monomial vector $y_i\cd e_\ell$ of $H_i\cd e_\ell$ by $\ol B_{i\ell}$ produces in the coefficient of $x_n^j\cd e_m$ the summand $u_{i\ell mj}$. The order $<_\zeta$ has to be chosen so that this monomial is the smallest one among the monomials of this coefficient (after having set $x=0$ and $y=0$). \med

(ii) A general monomial vector $x^\rho y^\sigma\cd e_\ell$ of $H_i\cd e_\ell$ is a multiple of the leading monomial vector $y_\i\cd e_\l$ of $B_{\i\l}$ with scope $q+n+\i$ and contributes to the coefficient of $x_n^j\cd e_m$ (for some $1\leq m\leq r$ and $0\leq j\leq d_m-1$, and after having set $x=0$ and $y=0$) if and only if $\l=\ell$, $\rho=(0\to 0, \rho_n)$ with $\rho_n\leq j$ and $\sigma = e_\i$, say $x^\rho y^\sigma\cd e_\l=x_n^{\rho_n}y_\i\cd e_\l$. The only contributions can be constant multiples of $u_{\i\l m j'}$ with $j'+\rho_n=j$. Note then that for this to happen we must have $x_n^{\rho_n}y_\i\cd e_\l>_\e\inin(H_i\cd e_\l)$ (otherwise this monomial does not appear in $H_i\cd e_\l$). Therefore $<_\zeta$ should satisfy\med

\hs 2cm $u_{\i\l m j'}>_\zeta u_{i\l mj}$
\hfill
for  $j' \leq j$ and $x_n^{\rho_n}y_\i\cd e_\l>_\e y_i\cd e_\l$,\med

\hfill say $j' \leq j$ and  $x_n^j\cd \inin(H_\i\cd e_\l)>_\e x_n^{j'}\cd \inin(H_i\cd e_\l)$.\med

(iii) A general monomial vector $x^\rho y^\sigma\cd e_\ell$ of $H_i\cd e_\ell$ is a multiple of the leading monomial vector $x_n^{d_\k}\cd e_\k$ of $B_{\k}$ with scope $q+n $ and contributes to the coefficient of $x_n^j\cd e_m$ (after having set $x=0$ and $y=0$)  if and only if $\k= \ell$, $\rho=(0\to 0, \rho_n)$ with $\rho_n= d_\k +t$ for some $t\geq 0$ and  $\sigma =(0\to 0)$, say $x^\rho y^\sigma\cd e_\l=x_n^{\rho_n} \cd e_\k$. The only contributions can be constant multiples of $u_{\k m j'}$ with $t+j'= j$. Note then that we must have $x_n^{\rho_n} \cd e_\k>_\e\inin(H_i\cd e_\k)=y_i\cdot e_\k$ and therefore $<_\zeta$ should satisfy\med

\hs 2cm $u_{\k m j'}>_\zeta u_{i\k mj}$
\hfill
for  $j' \leq j$ and $x_n^{\rho_n} \cd e_\k>_\e y_i\cd e_\k$,\med

\hfill say $j' \leq j$ and $x_n^j\cd \inin(G_\k) >_\e x_n^{j'}\cd \inin(H_i\cd e_\k)$.\med

%-------------------------------------------------

{\it Linear terms of $\, V_{i\ell m}(0,0,u,v)$}: These occur after the first substitution step of the polynomial division as the coefficients of  $e_m$ (with $1\leq i\leq p$, $1\leq \ell\leq s$, $r+1\leq m\leq s$) when dividing $H_i\cd e_\ell$  by the vectors $B_{\i\l}$ and $B_\k$ with leading monomial vectors $y_\i\cd e_\l$ and $x_n^{d_\k}\cd e_\k$ and scopes $q+n+\i $, respectively $q+n $  (where $\i $, $\l$ and $\k$ vary in the ranges  $1\leq \i \leq p$, $1\leq \l \leq s$, $1\leq \k \leq r$). \med

Let $x^\rho y^\sigma\cd e_\ell$ be a monomial vector of the expansion of $H_i\cd e_\ell$, with  $\rho\in \N^n$, $\sigma\in \N^p$. If it is a multiple of the leading monomial vectors $y_\i\cd e_\l$, respectively $x_n^{d_\k}\cd e_\k$, of $B_{\i\l}$, respectively $B_\k$, subject to the correct scope conditions, it will be replaced in the  polynomial division by the according multiple of the tails $\ol B_{\i\l}$, respectively $\ol  B_\k$. After the substitution we have to look at the coefficient of $e_m$ and set $x=0$ and $y=0$. We distinguish three cases.\med

(i) The substitution of the monomial vector $y_i\cd e_\ell$ of $H_i\cd e_\ell$ by $\ol B_{i\ell}$ produces in the coefficient of $e_m$ the summand $v_{i\ell m}$. The order $<_\zeta$ has to be chosen so that this monomial is the smallest one among the monomials of this coefficient (after having set $x=0$ and $y=0$). \med

(ii) A general monomial vector $x^\rho y^\sigma\cd e_\ell$ of $H_i\cd e_\ell$ is a multiple of the leading monomial vector $y_\i\cd e_\l$ of $B_{\i\l}$ with scope $q+n+\i$ and contributes to the coefficient of $e_m$ (after having set $x=0$ and $y=0$) if and only if $\ell=\l$, $\rho=(0\to 0)$ and $\sigma = e_\i$, say $x^\rho y^\sigma\cd e_\l=y_\i\cd e_\l$. The only contributions can be constant multiples of $v_{\i\l m}$. For this to happen we must have $y_\i\cd e_\l>_\e\inin(H_i\cd e_\l)$ (otherwise this monomial does not appear in $H_i\cd e_\l$). Therefore $<_\zeta$ should satisfy\med

\hs 2cm $v_{\i\l m}>_\zeta v_{i\l m}$
\hfill
for  $y_\i\cd e_\l>_\e y_i\cd e_\l$,\med

\hfill say $\inin (H_\i\cd e_\l)>_\e\inin(H_i\cd e_\l)$.\med

(iii) A general monomial vector $x^\rho y^\sigma\cd e_\ell$ of $H_i\cd e_\ell$ is a multiple of the leading monomial vector $x_n^{d_\k}\cd e_\k$ of $B_{\k}$ with scope $q+n $ and contributes to the coefficient of $e_m$ (after having set $x=0$ and $y=0$)  if and only if $\k=\ell$, $\rho=(0\to 0, \rho_n)$ with $\rho_n= d_\k$ and  $\sigma =(0\to 0)$, say $x^\rho y^\sigma\cd e_\l=x_n^{d_\k} \cd e_\k$. The only contributions can be constant multiples of $v_{\k m}$. Note then that we must have $x_n^{d_\k} \cd e_\k>_\e\inin(H_i\cd e_\k)$ and therefore $<_\zeta$ should satisfy\med

\hs 2cm $v_{\k m}>_\zeta v_{i\k m}$
\hfill
for  $x_n^{d_\k} \cd e_\k>_\e\inin(H_i\cd e_\k)$,\med

\hfill say   $\inin(G_\k) >_\e \inin(H_i\cd e_\k)$.\med

%-------------------------------------------------

{\it Linear terms of $\, U_{kmj}(0,0,u,v)$}: These occur after the first substitution step of the polynomial division as the coefficients of  $x_n^j\cd e_m$ (with $1\leq k\leq r$, $1\leq m\leq r$, $0\leq j\leq d_m-1$) when dividing $G_k$ by the vectors $B_{\iota\lambda}$ and $B_\kappa$ with leading monomial vectors $y_\iota\cd e_\lambda$ and $x_n^{d_\kappa}\cd e_\kappa$ and scopes $q+n+\i$, respectively $q+n $  (where $\i $, $\l$ and $\k$ vary in the ranges  $1\leq \iota \leq p$, $1\leq \lambda \leq s$, $1\leq \kappa \leq r$). \med

Let $x^\rho y^\sigma\cd e_\l$ be a monomial vector of the expansion of $G_k$, with  $\rho\in \N^n$, $\sigma\in \N^p$. If it is a multiple of the leading monomial vectors $y_\i\cd e_\l$, respectively $x_n^{d_\k}\cd e_\k$, of $B_{\iota\lambda}$, respectively $B_\kappa$, subject to the correct scope conditions, it will be replaced in the  polynomial division by the according multiple of the tails $\ol B_{\iota\lambda}$, respectively $\ol  B_\kappa$. After the substitution we have to look at the coefficient of $x_n^j\cd e_m$ and set $x=0$ and $y=0$. We distinguish three cases.\med

(i) The substitution of the monomial vector $x_n^{d_k}\cd e_k$ of $G_k$ by $\ol B_k$ produces in the coefficient of $x_n^j\cd e_m$ the summand $u_{kmj}$. The order $<_\zeta$ has to be chosen so that this monomial is the smallest one among the monomials of this coefficient (after having set $x=0$ and $y=0$). \med

(ii) A general monomial vector $x^\rho y^\sigma\cd e_\k$ of $G_k$ is a multiple of the leading monomial vector $y_\i\cd e_\l$ of $B_{\i\l}$ with scope $q+n+\i$ and contributes to the coefficient of $x_n^j\cd e_m$ (after having set $x=0$ and $y=0$) if and only if  $\k=\l$, $\rho=(0\to 0,\rho_n)$ and $\sigma = e_\i$, say $x^\rho y^\sigma\cd e_\l=x_n^{\rho_n}y_\i\cd e_\l$. The only contributions can be constant multiples of $u_{\i\l mj'}$ with $\rho_n+j'=j$, say $\rho_n=j-j'$. For this to happen we must have $x_n^{\rho_n}y_\i\cd e_\l>_\e\inin(G_k)$ (otherwise this monomial does not appear in $G_k$). Therefore $<_\zeta$ should satisfy\med

\hs 2cm $u_{\i\l mj'}>_\zeta u_{kmj}$
\hfill
for  $j'\leq j$ and $x_n^{\rho_n}y_\i\cd e_\l>_\e x_n^{d_k}\cd e_k$,\med

\hfill say $j'\leq j$ and $x_n^j\cd\inin (H_\i\cd e_\l)>_\e x_n^{j'}\cd\inin(G_k)$.\med

(iii) A general monomial vector $x^\rho y^\sigma\cd e_\k$ of $G_k$ is a multiple of the leading monomial vector $x_n^{d_\k}\cd e_\k$ of $B_{\k}$ with scope $q+n$ and contributes to the coefficient of $x_n^j\cd e_m$ (after having set $x=0$ and $y=0$)  if and only if $\rho=(0\to 0, \rho_n)$ with $\rho_n\geq d_\k$ and  $\sigma =(0\to 0)$, say $x^\rho y^\sigma\cd e_\l=x_n^{\rho_n} \cd e_\k$ with $\rho_n=d_\k+t$ for some $t\geq 0$. The only contributions can be constant multiples of $u_{\k mj'}$ with $t+j'=j$. Note then that we must have $x_n^{d_\k+t} \cd e_\k>_\e\inin(G_k)$ and therefore $<_\zeta$ should satisfy\med

\hs 2cm $u_{\k mj'}>_\zeta u_{kmj}$
\hfill
for  $j'\leq j$ and $x_n^{d_\k+t} \cd e_\k>_\e\inin(G_k)$,\med

\hfill say $j'\leq j$ and  $x_n^j\cd \inin(G_\k) >_\e x_n^{j'}\cd \inin(G_k)$.\med

%-------------------------------------------------

{\it Linear terms of $\, V_{k m}(0,0,u,v)$}: These occur after the first substitution step of the polynomial division as the coefficients of $e_m$ (with $1\leq k\leq r$, $r+1\leq m\leq s$) when dividing $G_k$  by the vectors $B_{\i\l}$ and $B_\k$ with leading monomial vectors $y_\i\cd e_\l$ and $x_n^{d_\k}\cd e_\k$ and scopes $q+n+\i $, respectively $q+n $  (where $\i $, $\l$ and $\k$ vary in the ranges  $1\leq \i \leq p$, $1\leq \l \leq s$, $1\leq \k \leq r$). \med

Let $x^\rho y^\sigma\cd e_\k$ be a monomial vector of the expansion of $G_k$, with  $\rho\in \N^n$, $\sigma\in \N^p$. If it is a multiple of the leading monomial vectors $y_\i\cd e_\l$, respectively $x_n^{d_\k}\cd e_\k$, of $B_{\i\l}$, respectively $B_\k$, subject to the correct scope conditions, it will be replaced in the  polynomial division by the according multiple of the tails $\ol B_{\i\l}$, respectively $\ol  B_\k$. After the substitution we have to look at the coefficient of $e_m$ and set $x=0$ and $y=0$. We distinguish three cases.\med
 
(i) The substitution of the monomial vector $x_n^{d_k}\cd e_k$ of $G_k$ by $\ol B_k$ produces in the coefficient of $e_m$ the summand $v_{km}$. The order $<_\zeta$ has to be chosen so that this monomial is the smallest one among the monomials of this coefficient (after having set $x=0$ and $y=0$). \med

(ii) A general monomial vector $x^\rho y^\sigma\cd e_\k$ of $G_k$ is a multiple of the leading monomial vector $y_\i\cd e_\l$ of $B_{\i\l}$ with scope $q+n+\iota$ and contributes to the coefficient of $e_m$ (after having set $x=0$ and $y=0$) if and only if $\k=\l$, $\rho=(0\to 0)$ and $\sigma = e_\i$, say $x^\rho y^\sigma\cd e_\l=y_\i\cd e_\l$. The only contributions can be constant multiples of $v_{\i\l m}$. For this to happen we must have $y_\i\cd e_\l>_\e\inin(G_k)$ (otherwise this monomial does not appear in $G_k$). Therefore $<_\zeta$ should satisfy\med

\hs 2cm $v_{\i\l m}>_\zeta v_{k m}$
\hfill
for  $y_\i\cd e_\l>_\e \inin(G_k)$,\med

\hfill say $\inin (H_\i\cd e_\l)>_\e\inin(G_k)$.\med

(iii) A general monomial vector $x^\rho y^\sigma\cd e_\k$ of $G_k$ is a multiple of the leading monomial vector $x_n^{d_\k}\cd e_\k$ of $B_{\k}$ with scope $q+n $ and contributes to the coefficient of $e_m$ (after having set $x=0$ and $y=0$)  if and only if $\rho=(0\to 0, \rho_n)$ with $\rho_n= d_\k$ and  $\sigma =(0\to 0)$, say $x^\rho y^\sigma\cd e_\k=x_n^{d_\k} \cd e_\k$. The only contributions can be constant multiples of $v_{\k m}$. Note then that we must have $x_n^{d_\k} \cd e_\k>_\e\inin(G_k)$ and therefore $<_\zeta$ should satisfy\med

\hs 2cm $v_{\k m}>_\zeta v_{k m}$
\hfill
for  $x_n^{d_\k} \cd e_\k>_\e\inin(G_k)$,\med

\hfill say   $\inin(G_\k) >_\e \inin(G_k)$.\med

%-------------------------------------------------

This concludes the computation of the required inequalities for  the order $<_\zeta$ on $\N^q$.  It will be a monomial order on $\N^q$, where $q$ is the number of the variables $u$ and $v$, and has to be graded lexicographic subject to the following relations\med

%-------------------------------------------------

\hs 2cm $u_{\i\ell m j'}>_\zeta u_{i\ell mj}$
\hfill
if   $j' \leq j$ and  $x_n^j\cd \inin(H_\i\cd e_\ell)>_\e x_n^{j'}\cd \inin(H_i\cd e_\ell)$,\med

%-------------------------------------------------

\hs 2cm $u_{i\ell mj'}>_\zeta u_{kmj}$
\hfill
if  $j'\leq j$ and $x_n^j\cd\inin (H_i\cd e_\ell)>_\e x_n^{j'}\cd\inin(G_k)$,\med

%-------------------------------------------------

\hs 2cm $u_{ikmj'} <_\zeta u_{kmj}$
\hfill
if   $j' \geq j$ and $x_n^{j}\cd \inin(H_i\cd e_k)<_\e x_n^{j'}\cd \inin(G_k) $,\med

%-------------------------------------------------

\hs 2cm $u_{\k mj'}>_\zeta u_{kmj}$
\hfill
if  $j'\leq j$ and  $x_n^j\cd \inin(G_\k) >_\e x_n^{j'}\cd \inin(G_k)$,\med

%-------------------------------------------------

\hs 2cm $v_{\i\ell m}>_\zeta v_{i\ell m}$
\hfill
if  $\inin (H_\i\cd e_\ell)>_\e\inin(H_i\cd e_\ell)$,\med

%-------------------------------------------------

\hs 2cm $v_{i\ell m}>_\zeta v_{km}$
\hfill
if  $\inin (H_i\cd e_\ell)>_\e\inin(G_k)$,\med

%-------------------------------------------------

\hs 2cm $v_{ikm} <_\zeta v_{km}$
\hfill
if   $\inin(H_i\cd e_k)<_\e \inin(G_k) $,\med

%-------------------------------------------------

\hs 2cm $v_{\k m}>_\zeta v_{km}$
\hfill
if   $\inin(G_\k) >_\e \inin(G_k)$.\med

%-------------------------------------------------

The indices vary in the regions \med

\hs 2cm $1\leq i,\i\leq p$, \med

\hs 2cm $1\leq \ell\leq s$, \med

\hs 2cm $1\leq m\leq r$, \med

\hs 2cm $1\leq j, j'\leq d_m-1$ and \med

\hs 2cm $1\leq k,\k\leq r$ \med

for the $u$ variables, respectively in the regions\med

\hs 2cm $1\leq i,\i\leq p$, \med

\hs 2cm $1\leq \ell\leq s$, \med

\hs 2cm $r+1\leq m\leq s$ and \med

\hs 2cm $1\leq k,\k\leq r$ \med

for the $v$ variables. It is checked that the inequalities for $<_\zeta$ do not contradict each other, i.e., that there actually does exist a monomial order $<_\zeta$ fulfilling the eight conditions. \med

%-------------------------------------------------

We now extend $<_\eps$ to a monomial order $<_\xi$ on
$\N^{q+n+p}\times\{1\to s\}$ defined by \med

\hs 1cm $(\gamma,\a,\b,\ell)<_\xi
(\gamma',\a',\b',\ell')$ \hs.3cm if \hs.3cm  $(\abs\gamma,(\a,\b,\ell) 
,\gamma)<_{lex}(\abs{\gamma'},(\a',\b',\ell'),\gamma')$. \med

Here, $<_{lex}$ denotes the lexicographic order on $\N\times (\N^{n+p}\times
\{1\to s\})\times \N^q$, where $\abs\gamma$ and $\abs{\gamma'}$ are compared as
elements of $\N$ with the natural order, $(\a,\b,\ell)$ and $(\a',\b',\ell')$ as
elements of $\N^{n+p}\times\{1\to s\}$ with the order $<_\eps$, and $\gamma$
and $\gamma'$ as elements of $\N^q$ with respect to the order $<_\zeta$. 
%
%-------------------------------------------------
%         CHECK OF COMPATIBILITY OF RELATIONS (INCOMPLETE)
%-------------------------------------------------
%
\ignore

For illustration, take the second and the fifth equation for $\ell=\k$ and $i=\i$ replacing in the latter $j'$ by $j$ and $j$ by $j''$, say\med

%-------------------------------------------------

\hs 2cm $u_{\k m j'}>_\zeta u_{i\k mj}$
\hfill
for   $j' \leq j$ and $x_n^j\cd \inin(G_\k) >_\e x_n^{j'}\cd \inin(H_i\cd e_\k)$,\med

%-------------------------------------------------

\hs 2cm $u_{i\k mj}>_\zeta u_{\k mj''}$
\hfill
for  $j\leq j''$ and $x_n^{j''}\cd\inin (H_i\cd e_\k)>_\e x_n^{j}\cd\inin(G_\k)$.\med

%-------------------------------------------------

This gives $u_{\k m j'}>_\zeta u_{\k mj''}$ for $j'\leq j''$ and then ?\med

%-------------------------------------------------

\recognize
The inequalities which were imposed on $<_\zeta$ ensure that -- as shown above -- the initial monomials with respect to $<_\xi$ of the linear terms of $U_{i\ell mj}(0,0,u,v)$, $V_{i\ell m}(0,0,u,v)$, $U_{kmj}(0,0,u,v)$ and $V_{km}(0,0,u,v)$ are $u_{i\ell mj}$, $v_{i\ell m}$, $u_{kmj}$ and $v_{km}$. This was needed to show that $U$ and $V$ satisfy the properties of a mother code.\med 

%-------------------------------------------------

(g) We show that $u(x')$ and $v(x)$ are the baby series of $U$ and
$V$. By definition, $u(x')$ and $v(x)$ vanish at zero.  We have already seen in
part (d) above that $r_{i\ell}=R_{i\ell}(x,u(x'),v(x))$ and 
$r_k=R_k(x,u(x'),v(x))$ are zero. As  $u(x')$ does not depend on $x_n$ and $U$
does not depend on $v$ it follows from the decomposition of $\co(\widetilde I)$ that $U(x,u(x'))$ and $V(x,u(x'),v(x))$ are zero. This is what had to be
shown and concludes the proof of Theorem 10.1 in the $x_n$-regular case.\med

\big\goodbreak

%-----------------------------------------------
%           PROOF OF THM. 11.1 FOR x_n-REGULAR MODULES
%-----------------------------------------------

{\Bf 14. Proof of Theorem 11.1 for $x_n$-regular modules}\med

By Theorem 9.1 we may assume that the module $I$ is given by a minimal
standard basis $g_1\to g_r\in \kxx^s$ with initial monomial vectors
$x_n^{d_k}\cd e_k$. Let $(H,G)\in K[x,y]^p\times K[x,y]^{s\times r}$ be the
family code of $g_1\to g_r$ and let $h=(h_1\to h_p)$ be the baby series vector
of the mother code $H=(H_1\to H_p)\in K[x,y]^p$, so that $g_k=G_k(x,h(x))$.
\med

By Lemma 8.1 the submodule $\widetilde I=\,\<(y_i-h_i)\cdot e_\ell,g_k\>$ of
$K[[x,y]]^s$ equals $\<H_i\cdot e_\ell,G_k\>$. Let $<_\e$ be an extension of
$<_\eta$ to $\N^{n+p}\times\{1\to s\}$ with $y_i \cd e_\ell<_\e\ x_j\cd
e_\ell$ for all $i$, $j$ and $\ell$ as defined in Lemma 8.2. By Theorem 10.1 in the $x_n$-regular case we may assume that we already dispose of a reduced standard basis $b_{i\ell},\,b_k$ of $\widetilde I$ with initial
monomial vectors $y_i\cdot e_\ell$ and $x_n^{d_k}\cdot e_k$ with respect to
$<_\eps$. The father code of $b_{i\ell},\,b_k$ is given by the virtual reduced
standard basis $B_{i\ell},\,B_k$ of $\widetilde I$, the mother code is the vector
$(U,V)$ of components $U_{i\ell mj}$, $V_{i\ell m}$, $U_{kmj}$ and $V_{km}$. We
denote by $(u(x'),v(x))$ with components $u_{i\ell mj}(x')$, $v_{i\ell m}(x)$,
$u_{kmj}(x')$ and $v_{km}(x)$ the corresponding baby series vector. \med

We wish to divide an algebraic power series vector $f\in\kxx^s$ by the
submodule $I=\,\<g_k\>$ of $\kxx^s$. We may assume that $f$ has the same baby series
vector $h$ as $g_1\to g_r$. Write $f=F(x,h(x))\in K[x,h]^s$ with father code
$F\in K[x,y]^s$. We divide $F$ by the  polynomial vectors $B_{i\ell}$ and $B_k$
according to the polynomial division algorithm (Theorem 4.4) with leading monomial vectors $y_i\cdot e_\ell$ and $x_n^{d_k}\cdot e_k$ and scopes $q+n+i$, respectively $q+n$ (we recall that $n$ is the number of $x$-variables, $q$ is the number of $u$- and $v$-variables). We get a decomposition\med

\hs 3cm$F=\sum \widetilde A_{i\ell}\cdot B_{i\ell} +\sum \widetilde A_k\cdot B_k +
C$\med

with some polynomials $\widetilde A_{i\ell}$ in $K[u,v,x,y]$, $\widetilde A_k$ in $K[u,v,x]$, and a polynomial vector $C\in K[u,v]\otimes\co(\widetilde I)$. Replacing in 
this equation $y$ by $h(x)$, $u$ by $u(x')$ and $v$ by $v(x)$ yields a
decomposition \med

\hs 3cm $f=\sum \tilde a_{i\ell}\cd \tilde b_{i\ell} +\sum \tilde a_k\cd b_k + c$\med
  
for some algebraic power series $\tilde a_{i\ell}, \tilde a_k\in \kxx$ and an algebraic power series vector $c\in \kxx^s$. The vectors $\tilde b_{i\ell}$ and $b_k$ are obtained from $B_{i\ell}$ and $B_k$ by substitution of the variables. \med

(a) The vector $c$ has mother code $H$, $U$ and $V$ and father code $C$. Expand $C$ into \med

\hs 1cm $C= \sum_{m=1}^r\sum_{j=0}^{d_m-1} C_{mj}(u,x')\cdot
x_n^j\cdot e_m + \sum_{m=r+1}^s C_{m}(u,v,x)\cdot e_m$,\med

with polynomials $C_{mj}(u,x')$ and $C_{m}(u,v,x)$. Observe that, similarly as in section 13, part (c), the polynomials $C_{mj}(u,x')$ will not depend on $v$.  Substituting in $C$ the variables $u$ and $v$ by $u(x')$ and $v(x)$ we obtain for $c$ the decomposition\med 

\hs .5cm $ c= \sum_{m=1}^r\sum_{j=0}^{d_m-1} 
C_{mj}(u(x'),x')\cdot x_n^j\cdot e_m + \sum_{m=r+1}^s 
C_{m}(u(x'),v(x),x)\cdot e_m$.\med

Therefore $c\in\co(I)$ as required. \med\goodbreak

(b) We will show that the vectors $\tilde b_{i\ell}$ belong to the module $\<b_k\>$, thus getting a decomposition \med

\hs 3cm $f= \sum a_k\cd b_k + c$\med

for some power series $a_k\in \kxx$. To this end, recall that $\<(y_i-h_i)\cd e_\ell, g_k\>\ =\ \<b_{i\ell}, b_k\>$  (as submodules of $K[[x,y]]^s$)  and that the vectors $b_k$ do not depend on $y_i$. Thus the replacement of $y_i$ by $h_i$ does not affect them and gives $\<b_k\>\ \subset\ \<g_k\>$. As the initial modules of these two modules are equal (being generated by $x_n^{d_k}\cd e_k$ for $1\leq k\leq r$), the Division Theorem for power formal series yields equality $\<g_k\>\ =\ \<b_k\>$. This shows that the $b_{i\ell}$ belong to the submodule $\<(y_i-h_i)\cd e_\ell, b_k\>$ of $K[[x,y]]^s$. Therefore, replacing $y_i$ by $h_i$ in $b_{i\ell}$ yields $\tilde b_{i\ell}\in\  \<b_k\>\ \subset K[[x]]^s$. \med

(c) We finally show that the power series $a_k\in \kxx$ are algebraic and that their codes can be computed algorithmically. For this we will express constructively the father codes $B_{i\ell}$ of $\tilde b_{i\ell}$ in terms of $B_k$ and $H_i\cd e_\ell$. \med

The problem which we have to solve here is the following: Assume given a submodule $J$ of $K[[x]]^s$ generated by polynomial vectors $P_1\to P_r$, and let $Q$ be a polynomial vector. We use Algorithm 1.7.6 of [GP] computing the polynomial weak normal form of a polynomial with respect to a polynomially generated ideal in a power series ring, together with the comment at the bottom of page 58. By definition of the polynomial weak normal form [GP, def.÷  1.6.5], we get the construction of a decomposition $S  Q= \sum W_k P_k + R$ with polynomials $S$, $W_k$ and $R$ such that $S(0)\neq 0$,  where $R$ equals the  remainder of the formal power series division of $S  Q$ by $P_1\to P_r$. In case that $Q$ already belongs to the ideal generated by $P_1\to P_r$ in the power series ring, this decomposition specializes to $Q=\sum \widetilde W_k P_k$ with rational coefficients $\widetilde W_k=W_k/S$ in the localization of the polynomial ring at $0$.\med

Apply this technique to the polynomial vectors $B_{i\ell}$ and the submodule  $J=\ \<B_k, H_i\cd e_\ell, U\cd e_\ell, V\cd e_\ell\>$ of $K[[x,y,u,v]]^s$ (with the obvious abbreviations for $U$ and $V$). By definition, $J$ is generated by polynomial vectors. We have to check that $B_{i\ell}\in J$. For this, recall that $\widetilde I =\,\<H_i\cdot e_\ell, G_k\>\, = \,\<b_{i\ell}, b_k\>$ as submodules of $K[[x,y]]^s$ and that $\tilde b_{i\ell}  \in\ \<b_k\>$ in $K[[x]]^s$. Then, by construction of $U$ and $V$, we get the equalities\med

\hs 1cm $\<B_{i\ell}, B_k, H_i\cd e_\ell, U\cd e_\ell, V\cd e_\ell\>\ =  \ \<\tilde b_{i\ell}, b_k, H_i\cd e_\ell, U\cd e_\ell, V\cd e_\ell\>  $ \med

\hs 5.8cm $=\ \<b_k, H_i\cd e_\ell, U\cd e_\ell, V\cd e_\ell\> $\med

\hs 5.8cm $=\ \<B_k, H_i\cd e_\ell, U\cd e_\ell, V\cd e_\ell\>$\med

\hs 5.8cm $= J $.\med

We conclude that $B_{i\ell}\in J$. This shows that we can write $B_{i\ell}$ as a linear combination of the $B_k$, $H_\iota\cd e_\lambda$, $U\cd e_\lambda$, $V\cd e_\lambda$ with constructible rational power series coefficients, say\med

\hs 2cm $B_{i\ell} =\sum_{i\ell k} W_{i\ell k} B_k$\hs 1cm modulo $H$, $U$ and $V$,\med

where $W_{i\ell k}\in K[[x,y,u,v]]$ are rational functions. Upon replacing $y_i$ by $h_i$, $u$ by $u(x')$ and $v$ by $v(x)$ only the $B_k$ will subsist (the evaluations of the other polynomial vectors $H_\iota\cd e_\lambda$, $U\cd e_\lambda$, $V\cd e_\lambda$ vanish). This shows that the $W_{i\ell k}$ are the father codes of the coefficients $w_{i\ell k}$ in the linear combinations $\tilde b_{i\ell}=\sum_{i\ell k} w_{i\ell k} b_k$ expressing $\tilde b_{i\ell}$ in terms of $b_k$. The mother codes are the components of the polynomial vectors $H$, $U$ and $V$.\med

By definition of $a_k$ in terms of $\tilde a_{i\ell}$ and $\tilde a_k$ it now follows that the series $a_k$ are algebraic and that their family codes can be constructed by a finite algorithm. This establishes Theorem 11.1 for $x_n$-regular
modules.\med

\big\goodbreak

%-----------------------------------------------
%         PROOFS OF THM 10.1  AND 11.1  
%-----------------------------------------------

{\Bf 15. Proofs of Theorems 10.1 and 11.1 in the general case}\med

The idea for proving both theorems in the general case is to split a given
minimal standard basis of $I$ into two groups specified by the variables
appearing in their initial monomial vectors. The first group consists of
generators whose initial monomial vectors are pure $x_n$-powers. The remaining
generators have initial monomial vectors which involve also some
other variable. \med

So let be given, by Theorem 9.1, vectors $g_1\to g_r$ which form a minimal standard
basis of $I$. Adding suitable monomial multiples of the $g_k$ we may assume
that  $g_1\to g_r$ form a minimal Janet basis of $I$ with scopes $n_1\to n_r$.  We order $g_1\to g_r$ and permute the components of
$\kxx^s$ so that, for some $1\leq t\leq r$, the vectors $g_1\to g_t$ are
$x_n$-regular with initial monomial vectors $x_n^{d_k}\cdot e_k$, and so that
the initial monomial vectors of the remaining $g_{t+1}\to g_r$ involve at least
one of the variables $x_1\to x_{n-1}$. It is easy to see that the scopes $n_{t+1}\to n_r$ of $g_{t+1}\to g_r$ are all $<n$. This implies that \med

\hs 3cm $I=\sum_{k=1}^t \kxx\cdot g_k + \sum_{k=t+1}^r K[[x']]\cdot g_k$.\med

Therefore no $g_{t+1}\to g_r$ need to be multiplied in the subsequent divisions
by $x_n$. \med

By Theorem 10.1 in the $x_n$-regular case we may assume that $g_1\to g_t$
form already a {\it reduced} standard basis of the submodule $I_0=\,\<g_1\to
g_t\>$ of $\kxx^s$. By Theorem 11.1 in the $x_n$-regular case we know
how to divide $g_{t+1}\to g_r$ by $g_1\to g_t$ through a finite algorithm for
the respective family codes. This allows us to assume that $g_{t+1}\to g_r$
belong to\med

\hs 1.5cm $M=\co(I_0)=\sum_{m=1}^t\sum_{j=0}^{d_m-1} K[[x']]\cdot x_n^j\cdot
e_m+\sum_{m=t+1}^{s}\kxx\cdot e_m$.\med

It follows from the box condition that the initial monomial vectors of
$g_{t+1}\to g_r$ have their non-zero entry in the first summand\med

\hs 3cm 
$M_1=\sum_{m=1}^t\sum_{j=0}^{d_m-1} K[[x']]\cdot x_n^j\cdot e_m$ \med

of $M$. Setting $I'=\sum_{k=t+1}^r K[[x']]\cdot g_k$ we have $I'\subset M$ and 
$\inin(I')\subset M_1$. The monomial order on $\N^n\times \{1\to s\}$ induces via the inclusion $M\subset K[[x]]^s$ in a natural way an ordering of the monomial vectors in $M$. \med

We may now apply induction on $n$ as follows.  \med

First notice that $\inin(I')$, as a submodule of the free finite
$K[[x']]$-module $M_1$, satisfies again Hironaka's box condition with respect to the induced ordering of the variables.  Secondly, no
division occurs in the second summand  $M_2=\sum_{m=t+1}^{s}\kxx\cdot e_m$
of $M$. Therefore, by induction on the number of variables and discarding the
(irrelevant) fact that the summand $M_2$ is not finitely generated as
$K[[x']]$-module, we may assume to know how to construct the {\it reduced}
standard basis of the $K[[x']]$-submodule $I'$ of $M$ by a finite algorithm
on the level of codes. Notice that this  basis, when considered as
vectors in $\kxx^s$, remains reduced with respect to $g_1\to g_t$ because
its elements belong to $M=\co(I_0)$.\med

So we may assume that $g_{t+1}\to g_r$ already form a reduced standard basis 
of $I'$.  By induction on $n$ we may apply the division algorithm of Theorem 11.1 to
$I'$ as a submodule of $M$. Thus we know how to divide effectively algebraic
power series vectors in $M$ by $I'$. \med

Apply this to the tails $\ol g_k=x_n^{d_k}\cdot e_k-g_k$ of $g_1\to g_t$. They
belong to $M$ since $g_1\to g_t$ are a reduced standard basis of $I_0$ and
$M=\co(I_0)$. We divide these $\ol g_k$ by $I'$. This
allows us to assume from the beginning that $g_1\to g_t$ are reduced with respect
to $g_{t+1}\to g_r$, i.e., that $\ol g_k \in \co(I')$ for $1\leq k\leq t$. As
$I'\subset M=\co(I_0)$, the new $g_1\to g_t$ form again a reduced standard
basis (the module they generate may be different from $I_0$, but its initial
module is the same). In total, we have found the reduced standard basis $g_1\to
g_r$ of $I$. This proves Theorem 10.1.\med 

As for Theorem 11.1, any algebraic power series vector $f\in 
\kxx^s$ we wish to divide  by $I=\,\<g_1\to g_r\>$ can first be divided by
$I_0=\,\<g_1\to g_t\>$ using Theorem 11.1 in the $x_n$-regular case. It thus yields a
remainder in $M=\co(I_0)$. Then, using induction on $n$ and that $I'$
satisfies the box condition in $M$, we may divide this remainder as vector in
$M$ by $I'$. The resulting remainder can be interpreted, via the
inclusion of $M$ in $\kxx^s$, as a vector in $\co(I)\subset \kxx^s$. It will
coincide with the remainder of the formal power series division of $f$ by $I$ in
$\kxx^s$. It does not matter here that the second summand
$\sum_{m=t+1}^{s}\kxx\cdot e_m$ of $M$ is not finitely generated as
$K[[x']]$-module, because no division occurs in the last $s-r$ components of
$f$. \med

This establishes the division algorithm for algebraic power series vectors $f$
in $\kxx^s$ by submodules $I$ with box condition. Theorem 11.1 is proven.\med

%\[Comment (14) from April 11 2011: Notice that for its construction in the preprint we use three clues. (a) The chosen monomial is the minimal one w.r.t.÷ the given monomial order (so as to give the power series division). (b) The chosen monomial can be reinterpreted as the maximal one in the virtual basis with respect to a new monomial order (so as to yield a finite division algorithm when dividing polynomials by the virtual reduced standard basis).  (c) The complementary space to the initial module is ``multiplicatively closed'' in the sense that the product of the tails of two elements of the actual reduced standard basis either lies again entirely in the complement, or entirely inside the initial module. This is necessary so as to have the virtual division coincide with the actual one. Don't forget this last point when trying to extend the division to initial monomials like $xy$, where the complement $K[[x]]+K[[y]]$ does not have this closedness property.]
\med\big\goodbreak

%-----------------------------------------------
%      APPENDIX
%-----------------------------------------------
 
%\cl{\Bf APPENDIX}\big

%-----------------------------------------------
%       EXAMPLE  OF PACO FROM FEB 5 2005
%-----------------------------------------------
 
%\recognize

{\Bf 16. Example }\med

In this section we show in a concrete situation how the algorithms of Theorem 10.1
and 11.1 work in practice (for more examples, see [Wa]).  We will consider an ideal in three variables generated by algebraic power series involving a single baby series. Our objective will be the computation of the codes of the reduced standard basis of the ideal. As it will turn out, the reduced standard basis will consist of polynomials, so that, at the end, there will be no mother codes needed and the father codes of the basis coincide with the elements of the basis. Nevertheless, the example is significant, since it is not at all clear how to construct the codes of the reduced standard basis without using the techniques developed in the paper.\med

The example is chosen so as to illustrate the various aspects of the algorithm (reduction, division, passage to vectors, induction on the number of variables). Some steps could also be performed directly using some ad hoc tricks due to the simplicity of some of the generators of the ideals and modules involved. This will be indicated correspondingly. Nevertherless, all portions of the algorithm will show off at least once.\med

As a general rule, each step in the computations below will be followed by a renaming of the involved objects so as to keep the presentation as systematic as possible. In the subsequent step, letters will always refer to this renamed object and not to the original object defined at earlier stages of the exposition.\med

The initial variables will be denoted $x$, $y$ and $z$, corresponding to $x_1$, $x_2$ and $x_3$ in the text, with this ordering. This will affect $x_n$-regularity, being here first $z$-regularity, then, later, $y$-regularity and finally $x$-regularity. Also, the involved polynomial divisions will use this ordering of the variables.\med

The additional auxiliary variables appearing in the mother codes will be denoted by $t_1$, $t_2$, $\ldots$ (instead of $y_1$, $y_2$, $\ldots$ as in the text). The respective baby series will be $h_1$, $h_2$, $\ldots$ \med

%-------------------------------------------------
%           START OF EXAMPLE
%-------------------------------------------------

We consider the ideal $I$ in  $K[[x,y,z]]$ generated by three power series $g_1$, $g_2$, $g_3$ given as\med

\hs 3cm $g_1= z^2+xyz + {1\over 4} xyz^2 + \ldots=$\med

\hs 3.4cm $ =z^2+xyh(z)$,\med

\hs 3cm $g_2=yz + x^2z+y^2z$,\med

\hs 3cm $g_3=y^2+xyz$.\med

Here,  \med

\hs 3cm $h(z)= 1-\sqrt{1-z}={1\over 2}z+{1\over 8} z^2+\ldots$\med

 is the only involved baby series. Its mother code $H$ is taken as \med

\hs 3 cm $H= 2t-  t^2+z$ \med

(so that $h=h(z)$ is the unique formal power series solution of $H(z,t)=0$ satisfying $h(0)=0$.) Later on, when other mother codes will appear, we shall set $t=t_1$, $h=h_1$ and $H=H_1$. The father codes of $g_1$, $g_2$, $g_3$ are\med

\hs 3cm $G_1= z^2+xyt$\med

\hs 3cm $G_2=yz + x^2z+y^2z$,\med

\hs 3cm $G_3=y^2+xyz$.\med

The last two $G_2$ and $G_3$ do not involve $t$ because $g_2$ and $g_3$ are polynomials and hence $G_2=g_2$ and $G_3=g_3$. For our purposes it will be sufficient to have just one generator which is a true series. \med

We wish to compute the family codes of the reduced standard basis of $I=\,\<g_1,g_2,g_3\>\,\subset K[[x,y,z]]$ with respect to a given monomial order on $\N^3$. We shall choose the graded lexicographical order $<_\eta$ on $\N^3$ with $x>y>z$. This yields the initial monomials \med

\hs 3cm $\inin(g_1)= z^2 $\med

\hs 3cm $\inin(g_2)=yz$,\med

\hs 3cm $\inin(g_3)=y^2$.\med

It will turn out these do not yet generate the initial ideal $\inin(I)$ of $I$. The missing monomial is $x^4z$, which is the initial monomial of the element \med

\hs 3cm $ g_4= x^4z-x^3yz^2+x^4yh(z)$\med

of $I$ with father code\med

\hs 3cm $G_4= x^4z-x^3yz^2+x^4yt$.\med

Actually, $g_1\to g_4$ form a standard basis of $I$ with respect to $<_\eta$. This basis is obviously not reduced.\med

%-------------------------------------------------
%          EXAMPLE:       OVERVIEW
%-------------------------------------------------

{\bf Overview:} For the convenience of the reader, let us list the various steps which will appear in the calculations (below, ``computation of ...'' will always mean ``computation of the code of ...''.)\med

{\bf Step 1:} Computation of a standard basis of $I$. In addition to $g_1$, $g_2$, $g_3$ we will get a fourth generator $g_4$ of $I$, the one from above.\med

{\bf Step 2:} Specification of all $x_n$-regular elements of this basis and computation of the reduced standard basis of the ideal $I_1$ generated by them. Here, $x_n$ is $z$; as only $g_1$ is $z$-regular, $I_1=\,\<g_1\>$ is principal and its reduced standard basis can be computed with the algorithm of [AMR, Thm.÷ 5.5] or, equivalently,  as described in Theorem 10.1 above in the $x_n$-regular case for principal ideals. The reduced standard basis of $I_1$ will again be denoted by $g_1$. Its tail belongs to $\co(I_1)\isom K[[x,y]]^2$, where $\co(I_1)=K[[x,y]]\oplus K[[x,y]]z$ denotes the canonical monomial direct complement of $I_1$ in $K[[x,y,z]]$ with respect to the chosen monomial order. \med

{\bf Step 3:} Reduction of $g_2$, $g_3$, $g_4$ by $I_1=\,\<g_1\>$. This is the division of $g_2$, $g_3$, $g_4$ by $g_1$ with the algorithm of [AMR, Thm.÷ 5.6] or, equivalently, as described in Theorem 11.1 above in the $x_n$-regular case for principal ideals, $x_n$ being here $z$. The reduced series will again be denoted by $g_2$, $g_3$, $g_4$. \med

{\bf Step 4:} Interpretation of $g_2$, $g_3$, $g_4$ as vectors in $\co(I_1)\isom K[[x,y]]^2 $ and computation of the reduced standard basis of the submodule $I_2=\,\<g_2,g_3,g_4\>$ of $K[[x,y]]^2$ generated by them. By Step 1, the vectors $g_2$, $g_3$, $g_4$  already form a  standard basis of $I_2$, so they need not be completed again. Step 4 consists of four substeps.\med

{\parindent .7 cm

\litem{} {\bf Substep 4A:} Specification of all $y$-regular elements among $g_2$, $g_3$, $g_4$ and computation of the reduced standard basis of the submodule $I_3$ of $K[[x,y]]^2$ generated by these as described in Theorem 10.1 for the $x_n$-regular case (only $g_2$ and $g_3$ will be $y$-regular, so that $I_3=\,\<g_2,g_3\>$.) The reduced standard basis of $I_3$ will again be denoted by $g_2$, $g_3$. Its tails belong to $\co(I_3)\isom K[[x]]^3$, where $\co(I_3)=(K[[x]]\oplus K[[x]]y)\times  K[[x]]$ denotes the canonical monomial direct complement of $I_3$ in $K[[x,y]]^2$ with respect to the chosen monomial order.\med

\litem{}{\bf Substep 4B:} Reduction of $g_4$ by $I_3=\,\<g_2,g_3\>$. This is the division of $g_4$ by $g_2$, $g_3$ in $K[[x,y]]^2$ as described in Theorem 11.1 above in the $x_n$-regular case, $x_n$ being now $y$. The reduced vector will again be denoted by $g_4$.\med

\litem{}{\bf Substep 4C:} Interpretation of $g_4$ as a vector in $\co(I_3)\isom K[[x]]^3 $ and computation of the reduced standard basis of the submodule $I_4$ of $K[[x]]^3$ generated by it  as described in Theorem 10.1 in the $x_n$-regular case, $x_n$ being here $x$.  The situation will be so simple that the reduced standard basis of $I_4$ can be read off directly without using Theorem 10.1. It will again be denoted by $g_4$.\med

\litem{}{\bf Substep 4D:} Reduction of $g_2$, $g_3$ by $I_4=\,\<g_4\>$. This is the division of the tails $\ol g_2$, $\ol g_3$ of $g_2$, $g_3$ by $g_4$ in $K[[x]]^3$ as described in Theorem 11.1 in the $x_n$-regular case, $x_n$ being here $x$. Again, the situation will be so simple that the reduction can be read off without using Theorem 11.1. The reduced vectors will again be denoted by $g_2$, $g_3$.\med

}

The reduced standard basis of $I_2$ obtained in step 4 is thus $g_2$, $g_3$, $g_4$.\med

{\bf Step 5:} Reduction of $g_1$ by $I_2=\,\<g_2,g_3,g_4\>$. This is the division of the tail $\ol g_1$ of $g_1$ by $g_2$, $g_3$, $g_4$ in $K[[x,y]]^2$ as described in Theorem 11.1 in the general case. This step consists of 2 substeps.\med

{\parindent .7 cm

\litem{} {\bf Substep 5A:} Reduction of $g_1$ by $I_3=\,\<g_2,g_3\>$. This is the division of  the tail $\ol g_1$ of $g_1$ by $g_2$, $g_3$ in $K[[x,y]]^2$ as described in Theorem 11.1 in the $x_n$-regular case, $x_n$ being here $y$. The reduced vector will again be denoted by $g_1$. Its tail belongs to $\co(I_3)\isom K[[x]]^3$.\med

\litem{} {\bf Substep 5B:} Reduction of $g_1$ by $I_4=\,\<g_4\>$. This is the division of the tail $\ol g_1$ of $g_1$ by $g_4$ in $K[[x]]^3$ as described in Theorem 11.1 in the $x_n$-regular case, $x_n$ being here $x$. The reduced vector will again be denoted by $g_1$.\med

}

{\bf Conclusion:} The vectors $g_1$, $g_2$, $g_3$, $g_4$ obtained after step 5 now have to be reinterpreted as power series in $K[[x,y,z]]$. By construction, they form
the reduced standard basis of the ideal $I$ we started with.\med

%-------------------------------------------------
%          EXAMPLE:   STEP 1
%-------------------------------------------------

{\bf Computations:} We start now with the explicit description of the various  stages of the construction of the reduced standard basis of the ideal $I$.\med

{\bf Step 1:} Computation of a minimal standard basis of $I$.\med

Let $\widetilde I= \,\<H, G_k\>\,=\,\<t-h,g_k\>$ be the ideal of $K[[x,y,z,t]]$ associated to $I$ as in Lemma 8.1 (here, no $e_\ell$'s appear, since we work with ideals instead of modules; the index $k$ varies between $1$ and $3$). We may apply Mora's tangent cone algorithm or Lazard's homogenization method. Let $u$ be a homogenizing variable, and denote by $ H^h$, $G_k^h$ the homogenized polynomials of $H$ and $G_k$ in $K[x,y,z,t,u]$.  \med

We  extend the monomial order $<_\eta$ on $\N^3$ first to an order $<_\epsilon$ on $\N^4$ (the set of exponents of series in $K[[x,y,z,t]]$) such that $\inin_\eps H=t$ and $\inin_\eps G_k= \inin_\eta g_k$, and than $<_\eps$ to an order $<_h$ on $\N^5$ (the set of exponents of series in $K[[x,y,z,t,u]]$) such that $\pi(\lm_h(H^h))=\inin_\eps (H)$ and $\pi(\lm_h (G_k^h)= \inin_\eps (G_k)$, where $\lm_h$ denotes the leading monomial of a polynomial with respect to $<_h$ and $\pi: \N^4\times \N\map \N^4$. 
\med

A polynomial Gr\"obner basis with respect to $<_h$ on $\N^5$ of the
ideal $J\subset K[x,y,z,t,u]$ generated by
$ H^h$, $G_1^h$, $G_2^h$ and $G_3^h$ is given by \med

\hs 3cm $ ut-{1\over 2}uz-{1\over 2}t^2, uy^2+zyx, uzy+zy^2+zx^2, uz^2+tyx$, \med 
\hs 3cm $zy^3-z^2yx+zyx^2$, $t^2y^2+2tzyx-z^2yx$, $z^2y^2-ty^2x+z^2x^2$, \med 
\hs 3cm $t^2zy+2tzy^2-ty^2x+2tzx^2$, $t^2z^2+2t^2yx-tzyx$, \med 
\hs 3cm $
z^3yx-tzyx^2+tyx^3$,$zy^2x^2-z^2x^3+zx^4$, $ty^3x-tzyx^2+tyx^3$, \med 
\hs 3cm $
z^3x^3-ty^2x^3$, $t^2yx^3+2tzx^4-z^2x^4$, $uzx^4-z^2yx^3+tyx^4$, \med 
\hs 3cm $
tz^2yx^3-{1\over 2}t^2zx^4+2tzx^5-z^2x^5+{1\over 2}tyx^5$, \med 
\hs 3cm $
t^3zx^4-4t^2zx^5-8tzyx^5+4z^2yx^5+4tzx^6-z^2x^6$.\med 

Now substitute $u$ by $1$ and $t$ by $h(z)$ to get a standard basis of $I$. It is given by $g_1$, $g_2$ and $g_3$ as above and the series $g_4$,  with\med

\hs 3cm  $g_4=x^4z-x^3yz^2+x^4yh(z)=$\med 

\hs 3.4cm  $ =x^4z-x^3yz^2+x^4y({1\over 2}z+{1\over 8}z^2+\ldots)$ \med 

and initial monomial $x^4z$. This series has as father code the polynomial\med

\hs 3cm $G_4=x^4z-x^3yz^2+x^4yt$.\med

The standard basis shows that the ideal $I$ satisfies Hironaka's box condition with respect to a monomial order such that $x<y<z$. The initial ideal is generated by $z^2$, $yz$, $y^2$ and $x^4z$. Moreover, it can be seen that the series $g_1$, $g_2$, $g_3$, $g_4$ form a Janet basis of $I$ with scopes $3$, $2$, $2$ and $1$ respectively.\med

%-------------------------------------------------
%          EXAMPLE:   STEP 2
%-------------------------------------------------

{\bf Step 2:} Computation of the reduced standard basis of the ideal $I_1 =\,\<g_1\>$.\med

Clearly, $g_1=z^2+xyh$ is the only $z$-regular series among $g_1\to g_4$. We set $I_1 =\,\<g_1\>\, \subset K[[x,y,z]]$. The monomials $xyz^m$ appearing in $xyh(z)$ are multiples of the initial monomial $z^2$ of $g_1$, therefore $g_1$ is not reduced (or in Weierstrass form). Let us apply the algorithm described in Theorem 10.1 for
$x_n$-regular series in order to find a reduced standard basis of the
ideal $I_1$. This algorithm coincides with the algorithm in [AMR, Thm.÷ 5.5]. The minimal reduced  standard basis $b_{11}$, $b_1$ of the ideal
$\widetilde I_1 =\,\<H,G_1\>\,\subset K[[x,y,z,t]] $ has the following form (with
the notation of the proof of Theorem 10.1).\med

\hs 3cm $b_{11}=t-b_{11}^\circ -u_{1110}(x')-u_{1111}(x')z$, \med

\hs 3cm $b_1=z^2-b_1^\circ-u_{110}(x')-u_{111}(x')z$,\med

where $b_{11}^\circ, b_1^\circ $ belong to $ K\oplus Kz$, the letter $x'$ stands for the variables $(x,y)$, and $u_{1110}(x')$, $u_{1111}(x')$, $u_{110}(x')$, $u_{111}(x')$ are power series vanishing at 0. To simplify let us write\med

\hs 3cm $b_{}=t-b_{}^\circ -u_{0}(x')-u_{1}(x')z $,\med

\hs 3cm $c=z^2-c^\circ-w_{0}(x')-w_{1}(x')z$.\med

We first compute $b^\circ $ and $c^\circ$ by setting $x$ and $y$ equal to $0$ in
the ideal $\tilde I_1$. We get the ideal \med

\hs 1.5cm $\< H(0,0,z,t), G_1(0,0,z,t)\>\, = \,\< H, z^2\>\, =\,\<
t-h,z^2\>\, \subset K[[t,z]]$. \med

From the mother code of $h(z)$ we can compute its Taylor expansion up to any given degree. In this case we have $ h={1\over 2}z+{1\over 8}z^2+ \cdots $. It follows that the
(minimal) reduced standard basis of the  ideal $\,\<t-h,z^2\>\,$ is $ t-{1\over 2}z$ and $z^2$. This implies that $b^\circ =
{1\over 2}z$ and $c^\circ =0$.\med

Next we have to find the family code for  the series $u_0(x,y)$, $u_1(x,y)$, $w_0(x,y)$, $w_1(x,y)$. We will divide -- using the polynomial division -- the polynomials $H$ and $G_1$ by the virtual reduced standard basis\med

\hs 3cm $B_{}=t-b_{}^\circ -u_{0}-u_{1}z = t-{1\over 2}z-u_{0}-u_{1}z$,\med

\hs 3cm $C=z^2-c^\circ-w_{0}-w_{1}z= z^2-w_0-w_{1}z$\med

of the ideal $\widetilde I_1$ with initial monomials $t$ and $z^2$, where $u_0,u_1,w_0,w_1$ are now just unknowns. The remainders $R,S$ of these divisions are \med

\hs 1cm $R=(-2u_1+u_0+2u_0u_1+{1\over 4}w_1+u_1w_1
+ u_1^2w_1)z -2u_0+u_0^2+{1\over 4}w_0+u_1w_0+u_1^2w_0$,\med

\hs 1cm $S=({1\over 2}xy+xyu_1+w_1)z+xyu_0+w_0$.\med

Let $U_1$, $U_2$, respectively $W_1$, $W_2$, be the coefficients of $z$ and $1$ in $R$ and $S$. It is easy to prove that they form a mother code  with baby series $u_0(x,y),u_1(x,y),w_0(x,y)$ and $w_1(x,y)$. In the present example  the solutions vanishing at $0$ of this mother code can be described in an equivalent and more explicit way as follows. From the four equations $U_1=U_2=W_1=W_2=0$ we get  
\med

\hs 3cm $u_0(x,y)=w_0(x,y)=0$,\med

\hs 3cm $w_1(x,y)=- {1\over 2}xy  -xy u_1(x,y)$,\med

\hs 3cm $u_1(x,y)=- {1\over 16}xy + {1 \over 16}x^2y^2 -
{67 \over 1024}x^3y^3 + O(x^4y^4)$,\med

where the last series is the unique solution vanishing at $0$ of the
equation\med

\hs 3cm $H_2(x,y,z,t_2)=8xyt_2^3+12xyt_2^2+16(1+xy)t_2+xy=0$\med

in a new variable $t_2$. In this way, $H_2$ becomes the
mother code of the algebraic series $u_1(x,y)$, its father code being
the polynomial $t_2$. The  father code of $w_1(x,y)$ is
$- {1\over 2}xy -xy t_2$.\med

The reduced standard basis of the ideal $I_1=\,\<g_1\>$ is given by substituting in the polynomial $C=z^2-w_0-w_1z$  the variables $w_0$ and $w_1$ by the series $w_0(x,y)=0$ and $w_1(x,y)= - {1\over 2}xy -xy u_1(x,y)$. We get the algebraic series $z^2+({1\over 2}xy+xyu_1(x,y))z$  with father code $C(0,- {1\over 2}xy -xyt_2,x,y,z)=z^2+({1\over 2}xy+xyt_2)z$. We denote this series in the sequel again by $g_1$, and call its father code $G_1$. The corresponding baby series $u_1(x,y)$ is denoted by $h_2(x,y)$ with mother code $H_2(x,y,z,t_2)$ from above. For later reference we collect the new data in a  table. \med

\hs 3cm $g_1=z^2+({1\over 2}xy+xyh_2(x,y))z$,\med

\hs 3cm $G_1= z^2+({1\over 2}xy+xyt_2)z$,\med

\hs 3cm $h_2(x,y)= - {1\over 16}xy + {1 \over 16}x^2y^2 -
{67 \over 1024}x^3y^3 + \ldots$,\med

\hs 3cm $H_2(x,y,z,t_2)=8xyt_2^3+12xyt_2^2+16(1+xy)t_2+xy$.\med

Note here that the original baby series $h=h_1$ has been eliminated.  \med

%-------------------------------------------------
%          EXAMPLE:   STEP 3
%-------------------------------------------------

{\bf Step 3:}  Reduction of $g_2$, $g_3$, $g_4$ by $I_1=\,\<g_1\>$.\med

We will apply the algorithm described in the proof of Theorem 11.1
for $x_n$-regular series to divide $g_2$, $g_3$, $g_4$ by $g_1$. It will be useful to add a new variable $t_3$ and define\med

\hs 3cm $H_3(x,y,z,t_1,t_2,t_3)=t_3+{1\over 2}xy+xyt_2$. \med

In this setting $(H_1,H_2,H_3)$ is the mother code of the baby series 
$(h_1,h_2,h_3)$ where $h_1=h(z)$ and $h_2=u_1(x,y)$ have been previously
defined and where $h_3$ equals $w_1(x,y)$ from above. It is clear from $\inin(I_1)= \,\<z^2\>$ that   $g_2=yz + x^2z+y^2z$ and $g_3=y^2+xyz$ are already reduced with respect to $I_1$. Let us reduce $g_4$. We shall use polynomial division. Let $\widetilde I_1=\,\<B, C\>$ be the ideal in $K[[x,y,z,t_1,t_2,t_3]]$ associated to $I_1$ as in Lemma 8.1 (it is checked that this is exactly the ideal of the lemma), with virtual reduced standard
basis\med

\hs 3cm  $B=t_1-{1\over 2}z-t_2z$,\med

\hs 3cm  $ C= z^2-t_3z$.\med

Dividing the father code $G_1$ of $g_1$ by $B$ and $C$ with initial monomials $t_1$ and $z^2$ we get \med

\hs 3cm $ G_4=x^4z-x^3yz^2+t_1x^4y=$\med

\hs 3.4cm $=x^4yB-x^3yC+D_4$,\med

where $D_4=({1\over 2}yx^4+yx^4t_2+x^4-yx^3t_3)z$. Let us replace $G_4$ by $D_4$ and call it again $G_4$. It is the father code of a new algebraic series, denoted again by $g_4$, and defined by $g_4=G_4(x,y,z,h_1,h_2,h_3)$.  We have \med

\hs 3cm $g_4= ({1\over 2} yx^4+yx^4h_2+x^4-yx^3h_3)z$.\med

The series $g_2,g_3,g_4$ are now reduced with respect to $I_1=\,\<g_1\>$. For later reference we collect the actual data in a table. \med

\hs 3cm $g_1=z^2+({1\over 2}xy+xyh_2(x,y))z$,\med

\hs 3cm $G_1= z^2+({1\over 2}xy+xyt_2)z$,\med

\hs 3cm $g_2=G_2=yz + x^2z+y^2z$,\med

\hs 3cm $g_3=G_3= y^2+xyz$,\med

\hs 3cm $g_4= ({1\over 2} yx^4+yx^4h_2+x^4-yx^3h_3)z$,\med

\hs 3cm $G_4=({1\over 2}yx^4+yx^4t_2+x^4-yx^3t_3)z$,\med

\hs 3cm $h_1={1\over 2}z+{1\over 8} z^2+\ldots$,\med

\hs 3cm $h_2 = - {1\over 16}xy + {1 \over 16}x^2y^2 -
{67 \over 1024}x^3y^3 + \ldots$,\med

\hs 3cm $h_3 = - {1\over 2}xy -xy h_2 $,\med

\hs 3 cm $H_1= 2t_1-   t_1^2+z$, \med

\hs 3cm $H_2 =8xyt_2^3+12xyt_2^2+16(1+xy)t_2+xy$,\med

\hs 3cm $H_3=t_3+{1\over 2}xy+xyt_2$. \med

%-------------------------------------------------
%          EXAMPLE:   STEP 4
%-------------------------------------------------

{\bf Step 4:} Computation of the reduced standard basis of the submodule $I_2=\,\<g_2,g_3,g_4\>$ of $\co(I_1)\isom K[[x,y]]^2$. \med

The canonical direct monomial complement $\co(I_2)$ equals $K[[x,y]]\oplus K[[x,y]]z$ and is therefore isomorphic to $ K[[x,y]]^2$  as $K[[x,y]]$-module. The three series $g_2$, $g_3$, $g_4$ are mapped under this isomorphism onto the vectors\med

\hs 3cm $ g_2=(0, y+ x^2+y^2)$,\med

\hs 3cm $ g_3=(y^2, xy)$,\med

\hs 3cm $ g_4=(0, x^4+ {1\over 2} x^4y+x^4yh_2-x^3yh_3)$.\med

The monomial order $<_\eta$ on $\N^3$ induces via the inclusion $K[[x,y]]\oplus K[[x,y]]z \subset K[[x,y,z]]$ a monomial order, also denoted by $<_\eta$, on $\N^2\times \{1,2\}$. The respective initial monomial vectors are \med

\hs 3cm $ \inin(g_2)=(0, y )$,\med

\hs 3cm $ \inin(g_3)=(y^2, 0)$,\med

\hs 3cm $ \inin(g_4)=(0, x^4)$.\med

We see that $g_2$ and $g_3$ are $y$-regular, whereas $g_3$ is not. By the proof of  Theorem 10.1 we first treat the submodule generated by $g_2$ and $g_3$. \med

%-------------------------------------------------
%          EXAMPLE:   SUBSTEP 4A
%-------------------------------------------------

{\bf Substep 4A:} Computing the reduced standard basis of the submodule
$I_3= \< g_2,g_3 \>$ of $  K[[x,y]]^2$. \med

The vectors $g_2,g_3$ are not the reduced standard basis of $I_3$ but form
at least a minimal standard basis. The father codes of $g_2$ and $g_3$ are $G_2=(0,y+x^2+y^2)$ and $
G_3=(y^2,xy)$ respectively. They do not depend on the variables
$t_i$. From the proofs of Thms.÷ 2 and 3 follows that we have to
consider the virtual reduced standard basis of $\widetilde I_3=I_3$. Said differently, 
we do not need to consider the vectors $B_{i\ell}$. Thus \med

\hs 3cm $B_2=y^2\cdot e_1-b_2^\circ-\sum_{m=1}^2\sum_{j=0}^{d_m-1} u_{2mj}y^je_m$,\med

\hs 3cm $B_3=y\cdot e_2-b_3^\circ-\sum_{m=1}^2\sum_{j=0}^{d_m-1} u_{3mj}y^je_m$,\med

where $d_1=2,d_2=1$ and the vectors $b_2^\circ$,  $b_3^\circ$ belong to
$(K\times K) \oplus (Ky\times (0))$. The vectors $b_2^\circ$,  $b_3^\circ$ are obtained by specializing $x$ to $0$
in $G_2$ and $G_3$. From $G_2(0,y)=(y^2,0),
G_3(0,y)=(0,y+y^2)$ we conclude that $b_2^\circ=b_3^\circ=(0,0)$.\med

We then apply the polynomial division to reduce $G_2$ and $G_3$ by
the virtual reduced standard basis $B_2$ and $B_3$ of $I_3$ with initial monomial vectors $y^2\cd e_1$ and $y\cd e_2$. The corresponding remainders
are \med\hs 3cm $  ( (u_{311}u_{211} +u_{211}u_{220}+u_{211} + u_{210})y+
u_{210}u_{220} + u_{210} + u_{310}u_{211}, $ \med

\hs 9cm $ u_{320}u_{211}+u_{220}^2+u_{220}+x^2) $, \med

\hs 3cm $ ((u_{211}x+u_{311})y+u_{210}x+u_{310},
u_{220}x+u_{320})$.\med

Therefore, the system \med

\hs 3cm $u_{311}u_{211}
+u_{211}u_{220}+u_{211} + u_{210}=0$,\med

\hs 3cm $u_{210}u_{220} + u_{210}
+ u_{310}u_{211}=0$,\med

\hs 3cm $u_{320}u_{211}+u_{220}^2+u_{220}+x^2=0$,\med

\hs 3cm $u_{211}x+u_{311}=0$,\med

\hs 3cm $u_{210}x+u_{310}=0$,\med

\hs 3cm $u_{220}x+u_{320}=0$\med

is the mother code for the series $u_{210}(x)$, $u_{211}(x)$, $u_{220}(x)$,  $u_{310}(x)$, $u_{311}(x)$, $u_{320}(x)$. From this system we get \med

\hs 3cm $u_{210}(x) = u_{211}(x) =
u_{310}(x)=u_{311}(x) = 0$,\med

\hs 3cm $ u_{220} (x)= h_4(x)$, \med

\hs 3cm $u_{320}(x) = -h_4(x)x$,\med

where \med
 
\hs 3cm $h_4(x)=- {1 \over 2}+ \sqrt{ {1 \over 4}-x^2}=-x^2-x^4-2x^6-5x^8+O(x^{10})$\med

is the unique solution vanishing at $0$ of the equation \med
 
\hs 3cm $H_4=t_4^2+t_4+x^2=0$. \med

The reduced standard basis of the submodule $I_3=\,\<g_2,g_3\>$ of $K[[x,y]]^2$ is \med

\hs 3cm $(0,y-h_4(x))$,\med

\hs 3cm $(y^2,xh_4(x))$.\med

We denote these vectors again by $g_2$ and $g_3$.  For later reference we collect the actual data in a table. \med

\hs 3cm $g_1=z^2+({1\over 2}xy+xyh_2(x,y))z$,\med

\hs 3cm $G_1= z^2+({1\over 2}xy+xyt_2)z$,\med

\hs 3cm $g_2=G_2=yz + x^2z+y^2z$,\med

\hs 3cm $g_3=G_3= y^2+xyz$,\med

\hs 3cm $g_4= ({1\over 2} yx^4+yx^4h_2+x^4-yx^3h_3)z$,\med

\hs 3cm $G_4=({1\over 2}yx^4+yx^4t_2+x^4-yx^3t_3)z$,\med

\hs 3cm $h_1={1\over 2}z+{1\over 8} z^2+\ldots$,\med

\hs 3cm $h_2 = - {1\over 16}xy + {1 \over 16}x^2y^2 -
{67 \over 1024}x^3y^3 + \ldots$,\med

\hs 3cm $h_3 = - {1\over 2}xy -xy h_2 $,\med

\hs 3cm $h_4 =-x^2-x^4-2x^6-5x^8+ \ldots$,\med

\hs 3 cm $H_1= 2t_1-   t_1^2+z$, \med

\hs 3cm $H_2 =8xyt_2^3+12xyt_2^2+16(1+xy)t_2+xy$,\med

\hs 3cm $H_3=t_3+{1\over 2}xy+xyt_2$, \med

\hs 3cm $H_4=t_4^2+t_4+x^2=0$. \med

%-------------------------------------------------
%          EXAMPLE:   SUBSTEP 4B
%-------------------------------------------------

{\bf Substep 4B:} Reduction of $g_4$ by the submodule $I_3=\,\<g_2,g_3\>$ of $K[[x,y]]^2$. \med

We reduce the vector $g_4=(0,{1\over 2}yx^4+yx^4h_2+x^4-yx^3h_3)$
by $I_3=\,\< g_2,g_3\>$. We point out
that it is not enough -- as the special shape of $g_2=(0,y-h_4(x))$ may suggest -- to replace $y$ by $h_4(x)$ in $g_4$ because  the power
series $h_2(x,y)$ and $h_3(x,y)$ depend on $x,y$.\med

The virtual  reduced standard basis $ B_{i\ell} $, $B_2$, $B_3$ with
$i=2,3,4$, $\ell=1,2 $, of the submodule $\widetilde I_3=
\<H_i\cd e_\ell,G_2,G_3\>$ of $K[[x,y,z,t_2,t_3,t_4]]$ as in Lemma 8.1 is  \med

\hs 3cm $B_{i\ell} = t_i \cdot e_\ell -b_{i\ell}^{\circ} - \sum_{m=1}^2
\sum_{j=0}^{d_m-1} u_{i\ell mj}\cdot y^j\cdot e_m$, \med

\hs 3cm $B_2=(0,y-h_4(x))$,\med

\hs 3cm $ B_3=(y^2,xh_4(x))$, \med\goodbreak

using here the computation we made in Substep 4A. To calculate the reduced standard basis of $\widetilde I_3$ we use polynomial division: We divide $H_{i\ell}$, $i=2,3,4$, $\ell=1,2 $, and $G_2$ and $G_3$ by $B_{\iota \lambda}$, $\iota=2,3,4$, $\lambda=1,2 $, and $B_2$, $B_3$  with leading monomial vectors
$t_\iota\cdot e_\lambda, \, y\cdot e_2, \, y^2\cdot e_1$
respectively. From the remainders of these  divisions we get -- by a rather tedious computation -- a system defining the mother code for the series $u_{i\ell
mj}(x)$. Another, more direct computation then shows that this system  can be transformed into an equivalent  system of form $H_5=H_6=0$ where\med

\hs 1cm $H_5=8xt_4t_5^3+16xt_4t_5^2+(16+16xt_4)t_5+xt_4 $, \med

\hs 1cm  $H_6=(16+12t_4t_5x+8t_4t_5^2x+16xt_4+{1\over 32}t_4^3x^3-{3\over 4}t_4^2x^2 -
{1 \over 2}t_4^2t_5x^2)t_6 +   $ \med

\hs 8cm  $+ x^3t_4+ {1 \over 512}t_4^3x^5- {3 \over 64}t_4^2x^4  $,\med
 
and where we have set $t_5=u_{2220}, \, t_6=u_{2120}$. The baby series vector of the mother code $(H_5,H_6)$ will be denoted by $(h_5,h_6)$.\med
 
Now we can apply polynomial division to reduce $G_4$ with respect to
the virtual reduced standard basis $ B_{i\ell}$, $B_2$, $B_3$ of $\widetilde I_3$. The division gives\med

\hs 3cm $G_4 = \sum \widetilde A_{i\ell} B_{i\ell} + \widetilde A_2 B_2 +
\widetilde A_3 B_3 + C_4$,\med
 
\hs 3cm $C_4=(0,((t_4+t_4^2)t_5+1+ {1 \over 2}t_4^2+ {1 \over 2}t_4)x^4)$, \med 

where $C_4$ is the father code of the reduction of $g_4$ by $g_2$, $ g_3$. We denote this reduction again by $g_4$. It is the vector obtained from $C_4$ by
substituting the variables $t_4$, $t_5$ by the power series $h_4$ and $h_5$. \med

%-------------------------------------------------
%          EXAMPLE:   SUBSTEP 4C
%-------------------------------------------------

{\bf Substep 4C:} Computation of the reduced standard basis of the submodule $I_4=\,\<g_4\>$ of $\co(I_3)\isom K[[x]]^3$.\med

By Substep 4B we have achieved that $g_4$ belongs to the canonical
direct monomial complement \med

\hs 3cm  $\co(I_3)=(K[[x]]\oplus K[[x]]y) \times K[[x]]$ \med

of $I_3$ in $K[[x,y]]^2$.  We will identify $\co(I_3)$ with $K[[x]]^3$ as $K[[x]]$-modules. Thus\med

\hs 3cm   $g_4=(0,0,((h_4+h_4^2)(h_5+ {1\over 2}) + 1)x^4)$. \med

The reduced standard basis of $I_4$ is  $(0,0,x^4)$ since
$h_4 (0)=0$ implies that $((h_4+h_4^2)(h_5+ {1\over 2}) + 1)$ is invertible in $K[[x]]$. Here we could also apply the
algorithm of Theorem 10.1 to compute the reduced standard basis of $I_4$. In this case
the computations are trivial because the base ring is the
principal ideal domain $K[[x]]$. We set again $g_4=(0,0,x^4)$ with father
code $G_4=(0,0,x^4)$.\med

%-------------------------------------------------
%          EXAMPLE:   SUBSTEP 4D
%-------------------------------------------------

{\bf Substep 4D:} Reduction of $g_2$, $g_3$ by $I_4=\,\<g_4\>$.\med

We apply the division algorithm of Theorem 11.1 in order to divide the tails $\ol g_2$ and $\ol g_3$ of $g_2$ and $g_3$ by $g_4$ (this is sufficient since $\inin(g_2)$ and $\inin(g_3)$ do not contribute to the remainders.) As $\ol g_2$, $\ol g_3$ and $g_4$ belong to $\co(I_3) =(K[[x]]\oplus K[[x]]y) \times K[[x]]$ we may treat them as vectors in $K[[x]]^3$.  We thus have\med

\hs 3cm  $g_2=(0,0,h_4)$,  \med

\hs 3cm  $g_3=  (0,0,-xh_4)$,\med

\hs 3cm  $g_4= (0,0,x^4)$.\med

Using that $h_4 = -x^2-x^4-2x^6 - \cdots$ it can be seen by inspection that the remainders of the division of $\ol g_2$ and $\ol g_3$ by $g_4$ are $(0,0, -x^2)$ and $(0,0, x^3)$. \med

This can also be seen alternatively by applying the polynomial division. Namely, as  
$\ol g_2$, $\ol g_3$ as well as $g_4$ belong to $(0)\times (0)
\times K[[x]]$ we can work with the respective last components in $K[[x]]$.
Let us consider the polynomials $\ol G_2=t_4$, $\ol G_3=-xt_4$ and $G_4=x^4$ as the father codes of the last components of $\ol g_2$, $\ol g_3$ and $g_4$ respectively.\med

Since the only baby series appearing in $g_2$, $g_3$, $g_4$ is $h_4 $ 
we have to consider the virtual reduced standard basis
$B_{41}, B_4$ of the ideal $\widetilde{I}_5 \subset K[[x,t_4]]$
generated by $H_4=t_4^2+t_4+x^2$ and $G_4=x^4$. One has \med

\hs 3cm $B_{41} = t_4-u_{4110}-u_{4111}x+u_{4112}x^2+u_{4113}x^3$,\med

\hs 3cm $B_4=x^4-u_{410}-u_{411}x+u_{412}x^2+u_{413}x^3$.\med

We get $u_{410}=u_{411}=u_{412}=u_{413}=0$ since the
initial monomial of the baby series $b_4$ of $B_4$ should be $x^4$. On the
other hand, the remainder of  the polynomial division of $H_4$ by
$B_{41}$ and $B_4$ is \med

\hs 3cm $(u_{4113} + 2 u_{4111} u_{4112})x^3 + (1 +
u_{4112} + u_{4111} )x^2 + u_{4111}x$,\med

which implies $u_{4111}=u_{4113}=0$ and $u_{4112}=-1$. The reduced standard
basis of $\widetilde{I}_5$ is $b_{41}=t_4+x^2$ and $b_4=x^4$. Finally, we have to divide $\ol G_2$ and $\ol G_3$ by $B_{41}$ and $B_4$ using the polynomial division with leading monomial vectors $t_4$ and $x^4$. One has\med

\hs 3cm  $\ol G_2=t_4=B_{41}-x^2$, \med

\hs 3cm $ \ol G_3=-xt_4=-x
B_{41}+ x^3$. \med

Rephrasing everything as vectors in $K[[x,y]]^2$,  the reductions of $g_2, g_3$ by $g_4$ are \med

\hs 3cm $(0,y+x^2)$,\med

\hs 3cm $(y^2,-x^3)$.\med

We set again $g_2=(0,y+x^2)$, $g_3=(y^2,-x^3)$, rewrite $g_4$ as $g_4=(0,x^4)$, together with their father codes
$G_2=(0,y+x^2)$, $ G_3=(y^2,-x^3)$ and $G_4=(0,x^4)$. This is the reduced standard basis of $I_3$; it coincides with what we have got at the beginning of this  substep.\med

{\bf Conclusion of Step 4:} To finish Step 4 we have to rewrite the preceding vectors as
algebraic power series in $x$, $y$, $z$ in order to obtain the reduced
standard basis of the ideal $I_2=\,\<g_2,g_3,g_4\>$ of $K[[x,y,z]]$. The corresponding reduced standard basis is given by the polynomials (we write again  $g_1$, $g_2$ and $g_3$)\med

\hs 3cm $g_2=yz+x^2z$,\med

\hs 3cm $g_3=y^2-x^3z$,\med

\hs 3cm $g_4=x^4z$. \med

They coincide with their father codes.\med

%-------------------------------------------------
%          EXAMPLE:   STEP 5
%-------------------------------------------------

{\bf Step 5:}  Reduction of $g_1$ by the submodule $I_2=\,\<g_2,g_3,g_4\>$ of $K[[x,y]]^2$.  \med

Recall that $g_1= z^2+({1\over 2}xy+xyh_2)z$. It suffices to divide the  tail $\ol g_1= -({1\over 2} xy+xyh_2)z$ of $g_1$  by $I_2=\,\<g_2,g_3,g_4\>$. We consider $\ol g_1 $
and $g_2$, $g_3$, $g_4$ as vectors in $\co(I_1)\isom K[[x,y]]^2$. Their father codes
are  $\ol G_1=(0,-xyt_2-{1\over 2}xy )$, $G_2=(0,y+x^2), \,
G_3=(y^2,-x^3),\, G_4=(0,x^4)$ respectively. The computation splits into two parts. Following Theorem 11.1 we will divide first
$\ol g _1$ by $I_3=\,\<g_2,g_3\>$ as vectors in $K[[x,y]]^2$ because $g_2$, $g_3$ are the $y$-regular power series among $g_2$, $g_3$, $g_4$. Afterwards, $\ol g _1$ will be divided by $I_4=\,\<g_4\>$ as a vector in $\co(I_3)\isom K[[x]]^3$. \med

Notice here that it is necessary to work with power series vectors in the canonical direct monomial complements $\co(I_1)$ and $\co(I_3)$. This is possible because, by the preceding steps,  $g_1$ is reduced with respect to itself (hence $\ol g_1$ belongs to $\co(I_1)$), $g_2$,  $g_3$ and $g_4$ are reduced with respect to $g_1$ (hence also belong to $\co(I_1)$), and $g_4$ is reduced with respect to $g_2$ and $g_3$ (hence belongs to $\co(I_3)\subset \co(I_1)$).\med

%-------------------------------------------------
%          EXAMPLE:   SUBSTEP 5A
%-------------------------------------------------

{\bf Substep 5A:}  Reduction of $g_1$ by $I_3=\,\<g_2,g_3\>$.  \med

We have to divide $\ol g_1$ by $g_2$ and $g_3$ as described in Theorem 11.1. We will use the polynomial division to divide the father code $\ol G_1$ of $\ol g_1$ by
the virtual reduced standard basis $ B_{21}$, $ B_{22}$, $B_2$, $B_3$  of $\widetilde{I}_3 =\,\< H_2\cdot e_1, H_2\cdot
e_2,  G_2,G_3\>$ in $ K[[x,y,t_2]]^2$. Notice
that   the only baby series appearing in $\ol g_1$,
$g_2$, $g_3$ is $h_2$. Therefore,  the only mother code appearing in $\widetilde{I}_1$
is $H_2$. We have \med

\hs 3cm $B_{21}=(t_2-u_{2110}-u_{2111}y,-u_{2120})$,\med

\hs 3cm $B_{22}=(-u_{2210}-u_{2211}y,t_2-u_{2220})$,\med

\hs 3cm $B_2=(0,y-u_{220})$,\med

\hs 3cm $B_3=(y^2,-u_{320})$,\med

where the form of $B_2$ and $B_3$ follows from the computation made in Substep 4D. The remainder of this polynomial division is $R=(0,x^3(u_{2220}+{1\over 2}))$. The algebraic series $u_{2220}(x)$ is defined by the mother code \med

\hs 3cm  $H_7=8x^3t_7^3+12x^3t_7^2-(16-16x^3)t_7+x^3$, \med

where we have set $t_7=u_{2220}$. This mother code $H_7$ results from the
division of $H_2\cdot e_\ell$ and $G_2$, $G_3$  by
$B_{i\lambda}$, $B_2$, $B_3$ and an appropriate simplification. Let us write $h_7 $ for the baby  series $u_{2220}(x)$ with mother code $H_7$. It then follows that the reduction of $\ol g_1=(0,-{1\over 2}xy-xyt_2)$
with respect to $I_3$  is $(0,x^3(h_7+{1\over 2}))$. We write this reduction again
as $\ol g_1=(0,x^3(h_7+{1\over 2}))$.  Note that it belongs to $\co(I_3)$.\med

%-------------------------------------------------
%          EXAMPLE:   SUBSTEP 5B
%-------------------------------------------------

{\bf Substep 5B:} Reduction of $g_1$ by $I_4=\,\<g_4\>$.\med

We have to divide the tail $\ol g_1$ of $g_1$ by $g_4$ as described in Theorem 11.1. For this we will consider $\ol g_1$ and $g_4$ as vectors in
$\co(I_3) \isom K[[x]]^3$. We have $\ol g_1 =(0,0,x^3(h_7+{1\over 2}))$
and $g_4=(0,0,x^4)$. Since $h_7(0)=0$ the reduction of
$\ol g_1$ with respect to $g_4$ is $(0,0,{1\over 2}x^3)$. As in
Substep 4D this reduction can be also computed by using the
polynomial division. We omit the details.\med

%-------------------------------------------------
%  EXAMPLE:    CONCLUSION
%-------------------------------------------------

{\bf Conclusion of example:}  Starting with the family code
$H_1=t_1^2-2t_1+z$, $G_1=z^2+xyt_1$, $G_2=yz+x^2z+y^2z$,
$G_3=y^2+xyz$ of the generators $g_1$, $g_2$ and $g_3$ of the ideal $I\subset K[[x,y,z]]$ with baby series $h=1-\sqrt{1-z}={1\over 2}z+{1\over 8} z^2+\ldots$ we have found the reduced standard basis of $I$ with respect to $<_\eta$ as the polynomials  (denoted again by $g_1$, $g_2$, $g_3$ and $g_4$) \med

\hs 3cm   $g_1=z^2-{1\over 2}x^3z$,\med

\hs 3cm  $g_2=yz+x^2z$,\med

\hs 3cm  $g_3=y^2-x^3z$,\med

\hs 3cm  $g_4=x^4z$.\med

They coincide with their father codes, and all baby series and mother codes have disappeared. We leave it as a challenge to the interested reader to find this basis of $I $ directly without using the algorithms of the paper.\med

\goodbreak

%-----------------------------------------------
%       EXTRAS
%-----------------------------------------------
 \ignore

{\bf Extras}\med

There remain several things which could or should be added.\med

(1) Different characterizations of the ring $A^h$ of algebraic power series, cf.
Iversen, Local genericity. Let $A=K[x]$ and $\hat A=K[[x]]$.\med

(a) $A^h$ is algebraic closure of $A$ in $\hat A$ (Nagata)?\med

(b) Azumaya's definition via idempotents.\med

(c) Grothendieck's definition as quotients $g=Q(h)/R(h)$ with $h=\ol z$ the
residue class of $z$ in $A[z]/P(z)$ for a polynomial $P$ with $P'(0)\neq
0$.\med

\med
%-----------------------------------------------

(2) Construction of mother code via normalization as in [AMR] and [BCR], cf.
the discussion at the Complu in Feb 2004 and the normalization algorithm of de
Jong and Greuel.\med

\med

%-----------------------------------------------

(3) Examples of the construction of the mother code of an algebraic series, of a reduced standard basis of a module and of the division algorithm via codes.\med

\med

%-----------------------------------------------

(4) Description how to determine from the mother code of an algebraic series
its initial monomial.\med

\med

%-----------------------------------------------

(5) Why not working with rational father codes? \med

\med

%-----------------------------------------------

(6) Is the composition of two algebraic series again algebraic? If yes,
what would be the code?\med

\med

%-----------------------------------------------

(7) Check if the various extensions of the monomial orders are ok.\med

\med

%-----------------------------------------------

(8) Determine the father codes $A_k$ of the quotients $a_k$ in the proof
of Theorem 11.1 in the $x_n$-regular case.\med

\med

%-----------------------------------------------

(9) Show that the quotients $a_k$ of the formal power series division
are algebraic series if they are subject to support conditions.\med

\med

%-----------------------------------------------

(10)  Possible better to omit the reference to Shiota. Is the
reference to Nash appropriate? Is Nagata the right reference for the Inverse
Function Theorem for algebraic series? The reference [CLO] p.÷ 222 is not correct\med

%----------------------------------------------- 
\recognize 

 \big\goodbreak

%-----------------------------------------------
%            REFERENCES
%-----------------------------------------------

{\Bf References} \med

{\parindent 1.1 cm

%\litem{[ACH]} Alonso, M.E., Castro-Jim\'enez, F., Hauser, H.: Polynomial echelons and babylonian division. Preprint.

\litem {[AHV]} Aroca, J.-M., Hironaka, H., Vicente, J.-L.: The theory of the 
maximal contact. Me\-morias Mat.÷ Inst.÷ Jorge Juan Madrid 29 (1975).

%\litem{[ALR]} [Delete?]  Alonso, M.E., Luengo, I., Raimondo, M.: An algorithm on quasi-ordinary polynomials. Proceedings AAECC 6. Lecture Notes in Comp.÷ Sci.÷ 357. Springer 1989.

\litem{[Am]} Amasaki, M.: Applications of the generalized Weierstrass preparation theorem to the study of homogeneous ideals. Trans.÷ Amer.÷ Math.÷ Soc.÷ 317 (1990), 1-43.

\litem{[AMR]} Alonso, M.E., Mora, T., Raimondo, M.: A computational model for algebraic  power  series. J.÷ Pure Appl.÷ Alg.÷ 77 (1992), 1-38.

\litem{[Ar1]} Artin, M.: Grothendieck topologies. Mimeographed notes, Harvard University  1962.

\litem{[Ar2]} Artin, M.: Algebraic approximations of structures over complete local
rings. Publ.÷ Math.÷ I.H.E.S.÷ 36 (1969), 23-58.

\litem{[AM]} Artin, M., Mazur, B.:  On periodic points. Ann.÷ Math.÷ 81 (1965), 82-99.

\litem {[BCR]} Bochnak, J., Coste, M., Roy, M.-F.: G\'eom\'etrie Alg\'ebrique
R\'eelle. Springer 1987.

\litem{[BG]} Baclawski, K., Garsia, A.: Combinatorial decomposition of a class of rings. Adv.÷ Math.÷ 39 (1981), 155-184.

\litem{[BK]} Bostan, A.,  Kauers, M.: The complete generating function for Gessel walks is algebraic. With an appendix by Mark van Hoeij. Proc.÷ Amer.÷ Math.÷ Soc.÷ 138 (2010), no. 9, 3063-3078.

\litem{[BM]} Bousquet-M\'elou, M., Mishna, M.: Walks with small steps in the quarter plane. Contemp.÷ Math.÷ 520 (2010), 1-39. 

\litem{[BP1]} Bousquet-M\'elou, M., Petkov\v sek, M.: Linear recurrences with constant coefficients. Discrete Math.÷ 225 (2000), 51-75.

\litem{[BP2]} Bousquet-M\'elou, M., Petkov\v sek, M.: Walks confined in a quadrant are not always $D$-finite. Theor.÷ Comp.÷ Sci.÷ 307 (2003), 257-276.

\litem{[BrK]} Brieskorn, E., Kn\"orrer, H.: Plane Algebriac Curves. Birkh\"auser 1986.

\litem{[CLO]} Cox, D., Little, J., O'Shea, D.: Using Algebraic Geometry. Springer
1998.

\litem{[dJP]} De Jong, T., Pfister, G.: Local Analytic Geometry. Springer 2000.

\litem{[Ga]} Galligo, A.: A propos du Th\'eor\`eme de Pr\'eparation de 
Weierstrass. Lecture Notes in Math.÷ 409, 543-579. Springer 1973.

%\litem{[Ga]} Galligo, A.: Th\'eor\`eme de division et stabilit\'e en g\'eom\'etrie analytique locale. Ann.÷ Inst.÷ Fourier  39 (1979), 107-184.

\litem{[GB]} Gerdt, V.,  Blinkov, Y.: Involutive bases of polynomial ideals. Simplification of systems of algebraic and differential equations with applications. Math.÷ Comput.÷ Simulation 45 no. 5-6 (1998), 519-541. 

\litem{[Gr1]} Gr\"abe, H.-G.: The tangent cone algorithm and homogenization. J.÷ Pure Appl.÷ Algebra 97 (1994), 303-312.

\litem{[Gr2]} Gr\"abe, H.-G.: Algorithms in local algebra. J.÷ Symb.÷ Comp.÷ 19 (1995), 545-557.

\litem {[Gra]} Grauert, H.: \"Uber die Deformation isolierter Singularit\"aten
analytischer Mengen. Invent.÷ Math.÷ 15 (1972), 171-198.

\litem{[GP]} Greuel, G.M., Pfister, G.: A Singular Introduction to Commutative Algebra. Springer, 2nd edition 2007.

\litem{[GPS]} Greuel, G.M., Pfister, G., Sch\"onemann, H.: Singular. A Computer Algebra System for Polynomial Computations. Universit\"at Kaiserslautern, www.singular.uni-kl.de.

\litem {[Hi1]} Hironaka, H.: Idealistic exponents of singularity. In: Algebraic
Geometry,  The Johns Hopkins Centennial Lectures. Johns Hopkins University Press
1977.

\litem {[Hi2]} Hironaka, H.:  Resolution of singularities of an algebraic variety over a field of characteristic zero.  Ann.÷  Math.÷ 79 (1964), 109-326.

\litem{[HM]} Hauser, H., M\"uller, G.: A rank theorem for analytic maps between
power series spaces. Publ.÷ Math.÷ I.H.E.S.÷ 80 
(1995), 95-115.

\litem {[Ja1]} Janet, M.: Sur les syst\`emes d'\'equations aux d\'eriv\'ees
partielles. J.÷ Math.÷ $8^{\rm e}$ s\'er., III (1920), 65-151.

\litem {[Ja2]} Janet, M.: Le\c cons sur les syst\`emes d'\'equations aux
d\'eriv\'ees partielles. Gauthiers-Villars 1929.

\litem {[KPR]} Kurke, H., Pfister, G., Rozcen, M.: Henselsche Ringe und algebraische Geometrie. Dt.÷ Verlag der Wissenschaften, Berlin 1975.

\litem {[Lf]} Lafon, J.-P.: S\'eries formelles alg\'ebriques. C.÷ R.÷ Acad.÷ Sci.÷
Paris 260 (1965), 3238-3241.

%\litem {[LM]} Lafon, J.-P., Marot, J.:   Alg\`ebre Locale. Collection Enseignement des Sciences. Hermann 2002.

\litem {[Lz]} Lazard, D.: Gr\"obner bases, Gaussian elimination and
resolutions of systems of algebraic equations. Computer algebra (London 1983). 
Lecture Notes Comp.÷ Sci.÷ 162, 146-156. Springer 1983.

\litem {[Mi]}  Mishna, M.: Classifying lattice walks restricted to the quarter plane. J.÷ Combin.÷ Theory Ser.÷ A 116 (2009), no. 2, 460-477.
 
\litem {[Mu]} Mumford, D.: The Red Book of Varietes and Schemes. Springer, 2nd edition 1999.

\litem {[Mo]} Mora, T.: An algorithm to compute the equations of tangent cones.
Proceedings Eurocam 82. Lecture Notes Comp.÷ Sci.÷ 144, 158-165. Springer 1982.

\litem {[Na1]} Nagata, M.: On the theory of henselian rings. Nagoya Math.÷ Journ.÷ 5 (1953), 5-57, and 7 (1954), 1-19.

\litem {[Na2]} Nagata, M.: Local rings. Interscience 1962.

\litem {[Ra]} Raynaud, M.: Anneaux locaux hens\'eliens.  Lecture Notes Math.÷ 169. Springer 1970.

\litem {[Re]} Rees, D.: A basis theorem for polynomial modules. Proc.÷ Cambridge Phil.÷ Soc.÷ 52 (1956), 12-16.

\litem {[Ri]} Riquier, C.: Les Syst\`emes d'\'Equations aux D\'eriv\'ees Partielles.
Gauthier-Villars 1910.

\litem {[Ru]} Ruiz, J.: The basic theory of power series. Vieweg 1993.

\litem {[Se1]} Seiler, W.: A combinatorial approach to involution and
$\delta$-regularity I: Involutive bases in polynomial algebras of solvable type.  Appl.÷ Algebra Engrg.÷ Comm.÷ Comput.÷ 20, no. 3-4 (2009), 207-259.

\litem {[Se2]} Seiler, W.: A combinatorial approach to involution and
$\delta$-regularity II: Structure analysis of polynomial modules with Pommaret bases.
Appl.÷ Algebra Engrg.÷ Comm.÷ Comput.÷ 20, no. 3-4 (2009), 261-338.

\litem {[Se3]} Seiler, W.: Spencer cohomology, differential equations, and Pommaret bases. In: Gr\"obner Bases in Symbolic Analysis, 169-216. Radon Ser.÷ Comput.÷ Appl.÷ Math.÷ 2. De Gruyter 2007. 

\litem {[SW]} Sturmfels, B., White, N.: Computing combinatorial decompositions of rings. Combinatorica 11 (1991), 275-293.

\litem {[Wa]} Wagner, D.: Algebraische Potenzreihen. Diploma Thesis, Univ.÷ Innsbruck 2005.

\litem {[Wi]} Wilczynski, E.J.: On the form of the power series for an algebraic function. Amer.÷ Math.÷ Monthly 26 (1919), 9-12.

\litem {[Za]}  Zariski, O.: Analytical irreducibility of normal varieties. Ann.÷ Math.÷ 49 (1948), 352-361. 

}

%-----------------------------------------------

\vskip 1.5cm

M.E.A.: Departamento de \'Algebra, \par
Universidad Complutense de Madrid,\par
mariemi@mat.ucm.es\med

F.J.C.J.: Departamento de \'Algebra, \par
Universidad de Sevilla,\par
castro@us.es\med

H.H.: Fakult\"at f\"ur Mathematik,  Universit\"at Wien, \par
Institut f\"ur Mathematik, Universit\"at Innsbruck,\par
herwig.hauser@univie.ac.at
\vfill\eject\end